\documentclass{amsart}
\usepackage{amssymb}
\usepackage{amsfonts}
\usepackage{latexsym}

\newtheorem{theorem}{Theorem}[section]
\newtheorem{lemma}[theorem]{Lemma}
\newtheorem{proposition}[theorem]{Proposition}
\newtheorem{corollary}[theorem]{Corollary} 
\theoremstyle{definition}  
\newtheorem{definition}[theorem]{Definition}
\newtheorem{example}[theorem]{Example}

\theoremstyle{remark}
\newtheorem{remark}[theorem]{Remark}

\newcommand{\id}{\mbox{id}}

\newcommand{\End}{\mbox{End}}

\newcommand{\Hom}{\mbox{Hom}} 
 
\renewcommand{\span}{\mbox{span}} 
    
\newcommand{\eps}{\varepsilon} 
\newcommand{\into}{\hookrightarrow}
\newcommand{\act}{\rightharpoonup} 

\newcommand{\lact}{\triangleright} 
 
\newcommand{\la}{\langle\,} 
\newcommand{\ra}{\,\rangle}

\newcommand{\1}{_{(1)}} 
\newcommand{\2}{_{(2)}} 
\newcommand{\3}{_{(3)}} 
\newcommand{\4}{_{(4)}}

\title[Frobenius extensions and weak Hopf algebras]
{Frobenius extensions and weak Hopf algebras}

\author{Lars Kadison}
\address{Matematiska Institutionen \\ G{\" o}teborg
University \\
S-412 96 G{\" o}teborg\\ Sweden}
\email{kadison@math.ntnu.no}
\thanks{The first author thanks A.A. Stolin and Uffe Haagerup for discussions,
and NorFA for support.}

\author{Dmitri Nikshych}
\address{U.C.L.A., Department of Mathematics, Los Angeles, CA 90095-1555, USA}
\curraddr{M.I.T., Department of Mathematics, Cambridge, MA 02139-4307, USA}
\email{nikshych@math.mit.edu}
\thanks{The second author was partially supported by the NSF grant DMS-9988796.
He is grateful to Leonid Vainerman for numerous valuable
discussions on weak Hopf algebras and thanks M.I.T.\ and Pavel Etingof  for
the kind hospitality during his visit}
 
\begin{document}

\begin{abstract}
We study a symmetric Markov extension of $k$-algebras
$N \into M$, a certain kind of Frobenius extension 
with conditional expectation that is tracial on the centralizer 
and dual bases with a separability property. 
We place a depth two condition on this extension, 
which is essentially the requirement that the Jones tower 
$N\into M \into M_1  \into M_2$ can be obtained by taking relative 
tensor products with centralizers $A = C_{M_1}(N)$ and $B = C_{M_2}(M)$.
Under this condition, we prove that  $N \into M$ is the invariant subalgebra 
pair of a weak Hopf algebra action by $A$, i.e., that $N = M^A$. 
The endomorphism algebra $M_1 = \End {}_N M$
is shown to be isomorphic to the smash product algebra $M \# A$. 
We also extend results of Szyma\'nski \cite{S}, Vainerman
and the second author \cite{NV1}, and the authors \cite{KN}.
\end{abstract}

\maketitle


\begin{section}
{Introduction and Preliminaries}
In its most general setting, the Jones tower is the iteration
of the endomorphism ring construction over
 any non-commutative ring extension $S \into R_0$, which results 
in a tower of rings over $R_0$ \cite{J2}.
The first step is to form $R_0 \into R_1 := \End_S^r R_0$ via 
the left regular representation. 
The process may then be repeated to obtain $R_1 \into R_2 := \End_{R_0} R_1$. 
For a  finite index subfactor \cite{J1} or  a Markov extension \cite{K2}
$N \subseteq M= M_0$ 
the algebras in the Jones tower  have their usual form 
$M_n = M_{n-1}e_n M_{n-1}$ for $n = 1,2,3,\ldots$ where $e_n$ are the Jones
idempotents.  Up to Morita equivalence of rings,
the Jones tower over a Markov extension has periodicity two.

Now {\em weak} Hopf algebras generalize Hopf algebras and   
 are Hopf-like objects with self-dual 
axioms, introduced by B\"ohm and Szlach\'anyi in \cite{BSz}
and in \cite{BNSz} with Nill.
It is well understood now that Hopf algebras and weak Hopf algebras
arise as non-commutative symmetries of Jones towers
of certain finite index inclusions of topological algebras over
the complex numbers.
For a finite index von Neumann subfactor $N \subseteq M$ it was shown
by Szyma\'nski \cite{S}  that the depth $2$ condition for the associated tower
of centralizers $C_M(N) \subseteq C_{M_1}(N) \subseteq C_{M_2}(M) 
\subseteq\cdots$ is equivalent to $A:=C_{M_1}(N)$ having a natural
structure of a Hopf $C^*$-algebra, if  $C_M(N) = \mathbb{C}1$. 
In the general case where $C_M(N)\supseteq \mathbb{C}1$, it was
shown by Vainerman and the second author 
\cite{NV1} that the depth two condition is equivalent to $A$
being a \textit{weak} Hopf $C^*$-algebra.  In both cases, 
$A$ acts on $M$ in such a way that $N = M^A$ and $M_1\cong M\# A$; 
moreover, $B:= C_{M_2}(M)$ is naturally identified with the (weak)
Hopf $C^*$-algebra dual to $A$. An outline of the proof of 
an analogous result for depth $2$ inclusions of unital $C^*$-algebras
was given very recently by Szlach\'anyi \cite{Sz}.

In \cite{NV2} it was shown 
that a  finite index and finite depth II$_1$ subfactor is embeddable 
in a weak Hopf $C^*$-algebra smash product inclusion; whence it is 
canonically determined via a Galois-type correspondence 
by some weak Hopf $C^*$-algebra  and its coideal $*$-subalgebra.
In this respect, weak Hopf $C^*$-algebras play the same role
as Ocneanu's paragroups \cite{Oc}.

In \cite{KN} hypotheses of  depth $2$ are placed on a
Markov extension $N \subseteq M $ of algebras over a field
$k$ with trivial centralizer
$C_M(N) = k1$ such that $A = C_{M_1}(N)$
can be given a semisimple Hopf algebra structure via the Szyma\'nski
pairing \cite{S}.  Moreover,
$A$ acts on $M$  in such a  way that the Jones tower above $M$ is isomorphic
to a duality-for-actions tower obtained from the smash product of $M$ and $A$
and the standard left action of $A^*$ on $A$:

\begin{equation}
\label{duality result}
\begin{array}{ccccccc}
N & \into & M  & \hookrightarrow &  M_1 & \hookrightarrow & M_2     \\
\|  &            &  \|  &        &  \downarrow {\scriptscriptstyle \cong} & 
          & \downarrow {\scriptscriptstyle \cong} \\ 
M^A   & \into & M & \into     &        M \# A &  \into & M \# A \# A^*. 
\end{array}
\end{equation}
We can continue iteration in the isomorphic copy of the Jones
tower  by alternately acting by $A$
and its dual $A^*$. Indeed,
it is a well-known theorem in algebra and operator algebras that
the  algebra $M \# A \# A^*$ above  is isomorphic to 
the endomorphism algebra $\End(M \# A)_M$ (cf.\
\cite{BM} for Hopf algebras and \cite{N}
for weak Hopf algebras). 

In this paper, we extend (\ref{duality result}) to
a Markov extension $N \into M$ which satisfies less restrictive conditions
than trivial centralizer and free extension $M_1/M$ as in \cite{KN}.
We assume conditions slightly stronger
than $U := C_M(N)$ being  a separable algebra on which the Markov trace
$T$ is non-degenerate. For the depth $2$ conditions, we assume
that the canonical conditional expectations $E_M$ and $E_{M_1}$
have dual bases in $A$ and its dual centralizer $B = C_{M_2}(M)$, 
respectively.
In exchange we obtain a canonical structure of a semisimple and cosemisimple 
weak Hopf algebra on $A$ with the dual $B$. Furthermore, the  
smash products above no longer have $k$-vector space structure given
by $M \# A = M \otimes_k A$ and $M \# A \# A^* \cong  M_1 \otimes_k B$, but
by $M \# A = M \otimes _U A$ and $M \# A \# A^* \cong M_1 \otimes_V B$,
where $V= C_{M_1}(M)$.
 
This paper is organized as follows.  

In this section we move on to cover preliminaries
essential to  this paper  --- weak Hopf algebras and their actions,
Markov extensions, the Basic Construction Theorem, and conditions of symmetry
and weak irreducibility on Markov extensions that  will be needed in the 
later sections.  

In Section~2  we place depth $2$ conditions on the Jones tower over 
a symmetric and weakly irreducible  Markov
extension, and develop a series of propositions and lemmas on depth $2$ properties
on the centralizers $U \subseteq A \subseteq C = C_{M_2}(N)$
and $V \subseteq B \subseteq C$, in both cases, $C$ 
being the basic construction
for Markov extensions of the same index as $M/N$. 

In Sections~3 and 4 we show that $A$ is a weak Hopf algebra
with the action on $M$ outlined  above.
First, in Section~3 we place a coalgebra
structure on $B$ by defining a non-degenerate pairing with $A$; the antipode
$S: B \to B$ comes from a symmetry in the definition of the pairing.
The rest of this section is devoted to proving that this structure on
$B$ satisfies the axioms of a weak Hopf algebra.  It follows that
$A$ is the dual weak Hopf algebra of $B$.  
Second, in Section~4 an action of $B$ on $M_1$
is introduced, and two equivalent expressions for this action are given.
Then we establish a left action of $A$ on $M$ with the outcome as in 
(\ref{duality result}):
the two vertical isomorphisms following from Theorems (\ref{phi theorem})
and (\ref{psi theorem}) together with Propositions (\ref{action of B on M_1})
and (\ref{A acts on M}), which establish the actions of $A$ and its dual. 

We note here that the main results in \cite[Sections 1-6]{KN} are recovered in
this paper if $U$ is trivial.  Furthermore, the results of this paper may 
be viewed as an answer to the question implicit 
in \cite[last line, p. 387]{BNSz}.  In an appendix, 
we extend to Markov extensions
the Pimsner-Popa formula for the Jones idempotent
generating the basic construction of composites in a Jones tower,
and also give a special example of a depth $2$ algebra extension.     

\subsection*{Weak Hopf algebras}

Throughout this paper we work over an arbitrary 
field $k$ and use the Sweedler notation for a 
comultiplication on a coalgebra $H$, writing $\Delta(h) = h\1 \otimes h\2$
for $h \in H$.

The following definition of a weak Hopf algebra and related notions were 
introduced in an equivalent form by B\"ohm and Szlach\'anyi in 
\cite{BSz} (see also \cite{BNSz}). We refer the reader to the recent 
survey \cite{NV3} for an introduction to the theory of weak Hopf algebras 
and its applications.

\begin{definition}[\cite{BSz}, \cite{BNSz}]
\label{weak Hopf algebra}
A {\em weak Hopf algebra} is  a $k$-vector space $H$ that has structures 
of an algebra $(H,\,m,\,1)$ and a coalgebra $(H,\,\Delta,\,\eps)$ 
such that the following axioms hold:
\begin{enumerate}
\item[1.] $\Delta$ is a (not necessarily unit-preserving) algebra
 homomorphism:
\begin{eqnarray}
\Delta(hg) &=& \Delta(h)\Delta(g).
\label{eqn: delta homo}
\end{eqnarray}
\item[2.] The unit and counit satisfy the identities:
\begin{eqnarray}
\eps(hgf) &=& \eps(hg\1)\eps(g\2f) =  \eps(hg\2)\eps(g\1f),
\label{eqn: eps m} \\
(\Delta \otimes \id) \Delta(1) &=&
(\Delta(1)\otimes 1)(1\otimes \Delta(1)) =   
(1\otimes \Delta(1))(\Delta(1)\otimes 1).
\label{eqn: delta 1}
\end{eqnarray}
\item[3.] There exists a linear map $S: H \to H$, called an {\em antipode}, 
satisfying the following axioms:
\begin{eqnarray}
m(\id \otimes S)\Delta(h) &=&(\eps\otimes\id)(\Delta(1)(h\otimes 1)),
\label{eqn: S epst} \\
m(S\otimes \id)\Delta(h) &=& (\id \otimes \eps)((1\otimes h)\Delta(1)),
\label{eqn: S epss} \\
S(h\1)h\2 S(h\3)&=& S(h),
\label{S id S}
\end{eqnarray}
\end{enumerate}
for all $h,g,f\in H$.
\end{definition}

Axioms (\ref{eqn: eps m}) and (\ref{eqn: delta 1}) 
are analogous to the bialgebra axioms specifying $\eps$ as an
algebra homomorphism and $\Delta$ a unit-preserving map, and 
Axioms (\ref{eqn: S epst}) and (\ref{eqn: S epss}) 
generalize the properties
of the antipode with respect to the counit. In addition, 
it may be shown
that given (\ref{eqn: delta homo}) - (\ref{eqn: S epss}),
Axiom (\ref{S id S}) is equivalent to $S$ being both an 
algebra and coalgebra anti-homomorphism.

A {\em morphism} of weak Hopf algebras is a map between them
which is both an algebra and a coalgebra morphism commuting with
the antipode. 

Below we summarize the basic properties of weak Hopf algebras,
see \cite{BNSz}, \cite{NV3} for the proofs.

The antipode $S$ of a weak Hopf algebra $H$ is unique; 
if $H$ is finite-dimensional then it is bijective.

The right-hand sides of the formulas (\ref{eqn: S epst}) 
and (\ref{eqn: S epss}) are called the {\em target} 
and {\em source counital maps} and denoted 
$\eps_t$, $\eps_s$
respectively:
\begin{eqnarray}
\label{counital maps}
\eps_t(h) = (\eps\otimes\id)(\Delta(1)(h\otimes 1)),\\
\eps_s(h) = (\id \otimes \eps)((1\otimes h)\Delta(1)).
\end{eqnarray}
The counital maps $\eps_t$ and $\eps_s$ are idempotents in $\End_k(H)$,
and satisfy  relations $S\circ \eps_t = \eps_s \circ S$ and
$S\circ \eps_s = \eps_t \circ S$.

The main difference between weak and usual Hopf algebras
is that the images of the counital maps are not necessarily equal 
to $k1_H$. They turn out to be  subalgebras of $H$ 
called {\em target} and {\em source counital subalgebras} 
or {\em bases} as they generalize the notion of a base
of a groupoid  (cf.\ examples below):
\begin{eqnarray}
H_t &=& \{h\in H \mid \eps_t(h) =h \}
  =  \{ (\phi\otimes \id)\Delta(1) \mid \phi\in H^* \}, \\
H_s &=& \{h\in H \mid \eps_s(h) =h \} 
  =  \{ (\id\otimes \phi)\Delta(1) \mid  \phi\in H^* \}.
\end{eqnarray}  
The counital subalgebras commute with each other 
and  the restriction of the antipode
gives an algebra anti-isomorphism between $H_t$ and $H_s$.


The algebra $H_t$ (resp.\ $H_s$) is separable (and, therefore,
semisimple) with the separability idempotent
$e_t = (S\otimes \id)\Delta(1)$ (resp.\ $e_s=(\id\otimes S)\Delta(1)$).

Note that $H$ is an ordinary Hopf algebra if and only if
$\Delta(1)=1\otimes 1$, iff $\eps$ is a homomorphism, 
and iff $H_t=H_s =k1_H$.

When $\dim_k H <\infty$,
the dual vector space $H^*=\Hom_k(H, k)$ has a natural structure of a weak
Hopf algebra with the structure operations dual to those of $H$:
\begin{eqnarray}
& & \la \phi\psi,\,h\ra = \la \phi\otimes\psi,\,\Delta(h) \ra, \\
& & \la \Delta(\phi),\,h \otimes g \ra =  \la \phi,\,hg \ra, \\
& & \la S(\phi),\,h\ra = \la \phi,\,S(h) \ra, 
\end{eqnarray}
for all $\phi,\psi \in H^*,\, h,g\in H$. The unit of $H^*$ is $\eps$
and the counit is $\phi \mapsto \la\phi,\, 1\ra$.


\begin{example}
\label{examples}
Let $G$ be a {\em groupoid} over a finite base (i.e., a category with 
finitely many objects, such that each morphism is invertible), 
then the groupoid algebra $kG$  is generated by morphisms $g\in G$ with 
the unit $1 =\sum_{X}\, \id_X$, where the sum is taken over all objects
$X$ of $G$, and  the product of  two morphisms is equal to  their composition
if the latter is defined and $0$ otherwise. 
It becomes 
a weak Hopf algebra via:
\begin{equation}
\Delta(g) = g\otimes g,\quad \eps(g) =1,\quad S(g)=g^{-1},\quad g\in G.
\end{equation}
The counital maps are given by $\eps_t(g) =gg^{-1}=\id_{{\rm target}(g)}$ and
$\eps_s(g) =g^{-1}g = \id_{{\rm source}(g)}$.

If $G$ is finite then 
the dual weak Hopf algebra $(kG)^*$ is generated by idempotents
$p_g,\, g\in G$ such that $p_g p_h= \delta_{g,h}p_g$ and
\begin{equation}
\Delta(p_g) =\sum_{uv=g}\,p_u\otimes p_v,\quad 
\eps(p_g)= \delta_{g,gg^{-1}}= \delta_{g,g^{-1}g},
\quad S(p_g) =p_{g^{-1}}.
\end{equation}
\end{example}
 
It is known that any group action on a set gives rise to a finite
groupoid. Similarly, in the non-commutative situation, one
can associate a weak Hopf algebra with every
action of a usual Hopf algebra on a separable algebra, see \cite{NTV}
for details. More interesting examples of weak Hopf algebras arise
from dynamical twisting of Hopf algebras \cite{EN},
closely related to the quantum dynamical Yang-Baxter equation,
and from the applications to the subfactor theory (\cite{NV1}, \cite{NV2}).

\begin{definition}[\cite{BNSz},  3.1]
\label{integral}
A left (right) {\em integral} in $H$ is an element
$l\in H$ ($r\in H$) such that 
\begin{equation}
hl =\eps_t(h)l, \qquad (rh = r\eps_s(h)) \qquad \mbox{ for all } h\in H. 
\end{equation}
 
\end{definition}
These notions clearly generalize the corresponding notions for Hopf 
algebras (\cite{M}, 2.1.1). We denote $\int_H^l$ (respectively, 
$\int_H^r$) the space of left (right) integrals in $H$ and by 
$\int_H = \int_H^l \cap \int_H^r$ the space of two-sided integrals.  

An integral in $H$ (left or right) is called {\em non-degenerate} if
it defines a non-degenerate functional on $ H^*$. A left integral $l$
is called {\em normalized} if $\eps_t(l)=1$. Similarly, $r\in \int_H^r$ 
is normalized if $\eps_s(r)=1$. The Maschke theorem for weak Hopf algebras
\cite{BNSz} states that a weak Hopf algebra $H$ is semisimple if and only
if it is separable, and if and only if it has a normalized integral. 
In particular, every semisimple weak Hopf algebra is finite dimensional.

\begin{example}
\label{examples of integrals}
(i) Let $G^0$ be the set of units of a finite groupoid $G$,
then the elements $l_e = \sum_{gg^{-1}=e}\,g\,(e\in G^0)$
span $\int_{kG}^l$ and elements $r_e = \sum_{g^{-1}g=e}\,g\,(e\in G^0)$
span $\int_{kG}^r$.
\newline
(ii) If $H = (kG)^*$ then $\int_H^l = \int_H^r =\span\{p_e| \, e\in G^0\}$.
\end{example}

\begin{definition}
\label{module algebra}
An algebra $A$ is a (left) {\em $H$-module algebra} if $A$ is a left
$H$-module via $h\otimes a \to h\cdot a$ and
\begin{enumerate}
\item[1)] $h\cdot ab = (h\1 \cdot a)(h\2 \cdot b)$,
\item[2)] $h\cdot 1 = \eps_t(h)\cdot 1$.
\end{enumerate} 
\end{definition}
If $A$ is an $H$-module algebra we will also say that $H$ acts on $A$.
The invariants $A^H = \{ a \in A |\, h \cdot a = \eps_t(h) \cdot a, 
\forall\, h \in H \}$ form a subalgebra by 2) above and a calculation 
involving  \cite[(2.8a),(2.7a)]{BNSz}. 

\begin{definition}
\label{comodule algebra}
An algebra $A$ is a (right) {\em $H$-comodule algebra} if $A$ is a right
$H$-comodule via $\rho: a \mapsto a^{(0)} \otimes a^{(1)}$ and
\begin{enumerate}
\item[1)] $\rho(ab) = a^{(0)} b^{(0)} \otimes a^{(1)} b^{(1)}$,
\item[2)] $\rho(1) = (\id\otimes\eps_t)\rho(1)$.
\end{enumerate} 
\end{definition}

The coinvariants $ A^{\text{co}H} = \{ a \in A |\, \rho(a) = a^{(0)} \otimes 
\eps_t(a^{(1)}) \}$ form a subalgebra of $A$.

It follows immediately that $A$ is a left $H$-module algebra
if and only if $A$ is a right $ H^*$-comodule algebra.

\begin{example}
\label{examples of actions}
\begin{enumerate}
\item[(i)] The target counital subalgebra
$H_t$ is a trivial $H$-module algebra via
$h\cdot z = \eps_t(hz)$, $h\in H,\,z\in H_t$.
\item[(ii)] $H$ is an $H^*$-module algebra via
the dual, or standard, action $\phi\act h = h\1 \la \phi,\,h\2\ra$, 
$\phi\in  H^*,\,h\in H$.
\item[(iii)] Let $A = C_H(H_s) =\{ a\in H \mid ay=ya \quad\forall y\in H_s\}$
be the centralizer of $H_s$ in $H$, then $A$ is an $H$-module algebra via
the adjoint action $ h\cdot a = h\1 a S(h\2)$.
\end{enumerate}
\end{example}

Let $A$ be an $H$-module algebra, then a {\em smash product}
algebra $A\# H$ is defined on a $k$-vector space $A\otimes_{H_t} H$,
where $H$ is a left $H_t$-module via
multiplication and $A$ is a right $H_t$-module via
$$
a\cdot z = S^{-1}(z)\cdot a = a(z\cdot 1), \qquad a\in A, z\in H_t,
$$
as follows. Let $a\# h$ be the class of $a\otimes h$ in $A\otimes_{H_t} H$,
then the multiplication in $A\#H$ is given by the familiar formula
$$
(a\# h)(b\# g) = a(h\1 \cdot b) \# h\2 g,\qquad a,b,\in A,\, h,g\in H,
$$
and the unit of $A\# H$ is $1\# 1$.

The smash product $A\# H$ is a left $H^*$-module algebra via
$$
\phi\cdot (a\#h) = a\#(\phi\act h), \quad
\phi\in H^*,\,h\in H,\,a\in A.
$$
It was shown in \cite{N} that there is a canonical isomorphism of algebras
$(A\#H)\#H^* \cong \End(A\#H)_A$, which extends the well-known duality
theorem for usual Hopf algebras \cite{BM}.

\subsection*{Symmetric Markov extensions}

Again let $k$ be a ground field. 
Recall that an algebra extension
$M/N$ is {\em Frobenius} if there is an $N$-bimodule homomorphism
 $E: M \rightarrow N$
and elements $\{ x_i \}$, $\{ y_i \}$ in $M$ such that for all $m \in M$, 
\begin{equation}
E(mx_i) y_i = m = x_iE(y_im),
\label{dual bases}
\end{equation}
where a summation over repeated indices is understood 
throughout the paper. We refer
to $E$, $\{ x_i \}$,
$\{ y_i \}$ as Frobenius coordinates,  $E$ being 
called a \textit{Frobenius homomorphism}, and the elements $\{ x_i \}$, $\{ y_i \}$ are
called \textit{dual bases}.  Another set of Frobenius coordinates 
$F: M \to N$,  $\{ r_j \}$, $\{ \ell_j \}$ is related to the first by
$F = Ed$ and \textit{ dual bases tensor} by
$e= r_j \otimes \ell_j = x_i \otimes d^{-1} y_i$
where $d = F(x_i) y_i$ is in the centralizer $C_M(N)$  \cite{Kasch1, O,W}.
Note that $e$ is a Casimir element, i.e., satisfies $me = em$ for all $m \in M$
by a computation as in Lemma (\ref{casimir}) below. 
A Frobenius homomorphism $E$ is left \textit{non-degenerate}
(or faithful) in the sense
that $E(xM) = 0$ implies $x = 0$; similarly, $E$ is right non-degenerate. 
Being Frobenius is a transitive property of extensions
with respect to the composition of Frobenius homomorphisms. 

An algebra extension $M'/N'$ is said to be \textit{split}
if $N'$ is isomorphic to an $N'$-bimodule direct summand in $M'$.  
For example, a Frobenius extension $M/N$ is split
if there is $d \in C_M(N)$ such that $E(d) = 1$
in the notation above, since $Ed$ is then a bimodule
projection $M \to N$.  

A Frobenius extension $M/N$ is \textit{symmetric} if there
is a Frobenius homomorphism $E$ such that $Eu = uE$
for each $u \in C_M(N)$; i.e., 
$E(ux) = E(xu)$ for all $x \in M, u \in C_M(N)$ \cite{Kasch2}. Let $U = C_M(N)$
for the rest of this section. For example,
 the symmetry condition is satisfied by a symmetric algebra $A/k$ 
\cite{Y}.
As an application of the symmetry condition,
we have: 

\begin{lemma}
\label{casimir}
For all $u \in U$,  
\begin{equation}
x_i u \otimes y_i = x_i \otimes u y_i
\end{equation}
in $M \otimes_N M$. 
\end{lemma}
\begin{proof}
We compute using Eqs.\ (\ref{dual bases}): 
\[
x_i u \otimes y_i = x_j E(y_j x_i u) \otimes y_i = x_j \otimes E(uy_j x_i) y_i
= x_j \otimes u y_j. \qed
\]
\renewcommand{\qed}{}\end{proof}  
 
Recall that a Frobenius extension $M/N$ is \textit{strongly separable}
if  $E(1) = 1$ and $x_iy_i = \lambda^{-1}1 \in k1$
for some $\lambda \in k^{\circ}$ whose reciprocal is called
the \textit{index}, denoted by $\lambda^{-1} = [M:N]_E$ \cite{K1, K2}. 
(In the terminology of \cite{W}, 
$E$ is a conditional expectation with quasi-basis $x_i, y_i$
such that the index of $E$ is nonzero in $k1$.)  
We say that a strongly separable extension is
a \textit{Markov extension}
if there is a (Markov) trace $T: N \rightarrow k$ such that $T(1) = 1_k$
and $T_0 := T \circ E$ is a trace on $M$ \cite{K1, K2}. 
A Frobenius homomorphism
$E$ that is a trace-preserving bimodule projection is referred to
as a \textit{conditional expectation}.

\begin{definition}
We refer to an extension of algebras $M/N$ as a \textit{symmetric Markov  
extension} if it is a Markov extension with coordinates
$E$, $\{x_i\}$, $\{y_i\}$ and Markov trace $T$ such that for each  $u \in U$,
 $E(ux) = E(xu)$ for every $x \in M$.
\end{definition}

 For example, the symmetric Frobenius condition  is  satisfied by the
irreducible Markov extensions in \cite{KN}, since $U$ is trivial for these.  
As another example, the symmetric Frobenius condition is satisfied  when $T$
is non-degenerate on $N$, e.g., for a II$_1$ subfactor
$N \subseteq M$ of finite index \cite{J1}: 

\begin{proposition}
If the Markov trace $T$ is non-degenerate on $N$, then
$uE = Eu$ for every $u \in U$.
\end{proposition}
\begin{proof}
We note that: for all $ n \in N, m \in M $
\[
  T(n E(um)) = T_0(num) = T_0(unm) = T_0(nmu) = T(nE(mu)), 
\]
which implies that $E(um) = E(mu)$ for all $m \in M$. 
\end{proof}

Let $M_1 =  M \otimes_N M \cong \End (M_N) $ denote the basic
construction of $M/N$: i.e., $M_1 = M e_1 M$ where $e_1 = 1 \otimes 1$ is 
the first
Jones idempotent  with conditional expectation
$E_M: M_1 \rightarrow M$ given by $E_M(me_1 m') = \lambda mm'$,
dual bases $\{ \lambda^{-1} x_i e_1 \},\ \{ e_1 y_i \}$, and index-reciprocal
$\lambda$. Recall that $M _1 \cong \End(M_N)  $ is given by 
$ m e_1 m'\mapsto \ell_m E \ell_{m'}$ where $\ell_m$ is left
multiplication by $m \in M$. The $E$-multiplication induced by composition
on $\End(M_N)$ is given by $$ e_1 m e_1 = e_1 E(m) = E(m) e_1 $$
for all $m \in M$. 

\begin{theorem}[``Basic Construction'']
Suppose $N \hookrightarrow M$ is a symmetric Markov extension.  
Then $M_1$ is a symmetric Markov extension of $M$ with Markov
trace $T_0 = T \circ E$ and is characterized by having idempotent
$e_1$  and conditional expectation $E_M: M_1 \to M$
such that 
\begin{enumerate}
\item $M_1 = M e_1 M$;
\item $ E_M(e_1) = \lambda 1$;
\item  for each $x \in M$: $e_1 x e_1 = e_1 E(x) = E(x) e_1$.
\end{enumerate}
\label{basic construction}
\end{theorem}
\begin{proof}
Most of the proof is found in \cite{K1} or \cite{K2}: we
 need only establish the symmetric Frobenius condition
 as well as the characterization above. 

Let $V = C_{M_1}(M)$.  Note that $U$ is anti-isomorphic to $V$
as algebras, via the map 
\begin{equation}
\label{phi}
\phi: U \rightarrow V,\ \ \phi(u) = x_i u e_1 y_i,
\end{equation}
which has inverse given by $v \mapsto \lambda^{-1} E_M(ve_1)$. Clearly
then $V \cong U^{\rm op}$.   
Note also that $$
E_M(ve_1) = E_M(e_1 v)$$ as a consequence of Lemma (\ref{casimir}).
 
We compute that $E_M v = v E_M$ for all $\phi(u) \in V$ and all $a,b \in M$:
\[
E_M(\phi(u) ae_1 b) = E_M(x_i u e_1 y_i a e_1 b) = E_M(x_i E(y_i a) u e_1 b)
= \lambda aub, \]
while
\[
E_M(ae_1b \phi(u)) =  E_M(ae_1 b x_i u e_1 y_i) = E_M(ae_1E(bx_i u) y_i)
= \lambda a u b. 
\]

Suppose $\tilde{M}$ is
an algebra with idempotent $f$ and conditional
expectation $\tilde{E}: \tilde{M} \to M$  satisfying
the conditions above. Since $\tilde{M} = M f M$ and $nf = fn$
for each $n \in N$, there is a surjective
mapping of $M_1 \to \tilde{M}$.  If 
$xf = fx$ for some $x \in M$, then $fxf= fE(x) = fx$; 
applying $\tilde{E}$ and the Condition (2),
 we see that $x = E(x) \in N$.  It follows
that the mapping $M_1 \to \tilde{M}$
is an  algebra isomorphism forming a commutative
triangle with  $\tilde{E}$ and $E_M$. 
\end{proof}

Let us recall that a  $k$-algebra $A$ is 
{\em Kanzaki  separable} (also
called strongly separable algebra in the literature) if
it has a symmetric separability element, or equivalently,
if the trace of the left regular representation of $A$ on itself
has dual bases $\{ x_i \}$ and $\{ y_i \}$ such that $x_i y_i = 1$ \cite{Kan}. 
Yet another equivalent condition:  $A$ is $k$-separable with invertible
Hattori-Stalling rank as a finitely generated projective module over its center \cite{SK}. 
For example, the full $p$-by-$p$ matrix algebra over a characteristic $p$ field $F$ 
is separable but not Kanzaki separable.  Over a non-perfect field $F$, a separable $F$-algebra 
is in turn  finite dimensional semisimple, but not necessarily the converse.  
In characteristic zero, all three notions coincide.   

For the rest of this paper, we will make the two assumptions
below on a symmetric Markov extension $M/N$.  
\begin{enumerate}
\item (Symmetric product assumption.)
$ x_iy_i =  y_i x_i = \lambda^{-1}1 \in k1$. 
\item (Weak irreducible assumption.)
$U$ is a Kanzaki separable $k$-algebra  with non-degenerate trace
$T_0|_U$.
\end{enumerate}

Under these assumptions,
it follows from the proof of the basic construction theorem that $V$ is also Kanzaki  
separable.  The next proposition shows
that $T_1 := T E E_M$ is a non-degenerate trace on $V$. 

\begin{lemma}
\label{V too}
We have the identity $T_1 \circ \phi = T_0$ on $U$. 
\end{lemma}
\begin{proof}
Let $u \in U$.  We compute using the symmetric product assumption:
\[
T_1 \phi(u) = T_1(x_i ue_1y_i) = \lambda T_0(x_i uy_i)=
\lambda T_0(y_i x_i u) =  T_0(u). \qed
\]

\renewcommand{\qed}{}\end{proof}  

\begin{remark}
For the purposes of this paper, 
the symmetric product assumption may be replaced by the identity
in the statement of Lemma (\ref{V too}).  
This last condition holds trivially
for an irreducible Markov extension as in \cite{KN}. 
\end{remark}

Since $M_1/M$ is also a symmetric
Markov extension with index $\lambda^{-1}$, 
we now iterate the basic construction to form $M_2 = M_1 e_2 M_1$
with conditional expectation $E_{M_1}(xe_2 y) = \lambda xy$
for each $x,y \in M_1$ and {\em second Jones idempotent}
$e_2$. We recall the
braid-like relations, $$e_1 e_2 e_1 = \lambda e_1$$ and $$e_2 e_1 e_2 = 
\lambda e_2$$  established in \cite{K2}, 
and the Pimsner-Popa relations, 
\begin{equation}
\label{PP}
xe_1 = \lambda^{-1} E_M(xe_1)e_1\ \ \ \ 
\forall \, x \in M_1,
\end{equation}
and three more similar equations for $e_1x$,
$e_2y$ and $ye_2$ where $y \in M_2$ \cite{KN}.

\end{section}

\begin{section}
{Properties of depth 2 extensions}

Let $M/N$ be a  symmetric Markov extension satisfying the weak irreducible
condition and the symmetric product condition in Section~1. 
Recall that this entails three conditions on a Markov extension
$(E: M \to N, x_i,y_i, \lambda, T: N \to k)$:
\begin{enumerate}
\item $E: M \to N$ is symmetric: $Eu = uE$ for each $u \in U = C_M(N)$.
\item $U$ is Kanzaki separable and $T_0|_U$ is a non-degenerate trace.
\item $ y_i x_i = \lambda^{-1} = x_i y_i$; alternatively, $T_0|_U = T_1 \circ 
\phi$ where $\phi: U \to V$ is the anti-isomorphism
defined in Eq.\ (\ref{phi}).
\end{enumerate}

In this section, we work with  the Jones tower above $M/N$: 
\begin{equation}
\label{Jones tower}
N \stackrel{\xleftarrow{E}}{\hookrightarrow} M \stackrel{\xleftarrow{E_M}}{\hookrightarrow} M_1 
\stackrel{\xleftarrow{E_{M_1}}}{\hookrightarrow} M_2.
\end{equation}
We denote the ``second'' centralizers by $A = C_{M_1}(N)$,
$B= C_{M_2}(M)$, and the ``big'' centralizer by 
$C= C_{M_2}(N)$, which contains $A$, $B$. Note that $U$ and $V$ 
are contained in $A$; 
$V$ and $W = C_{M_2}(M_1)$ are contained in $B$; and 
$V =  A \cap B$. See Figure 1. 

\begin{figure}
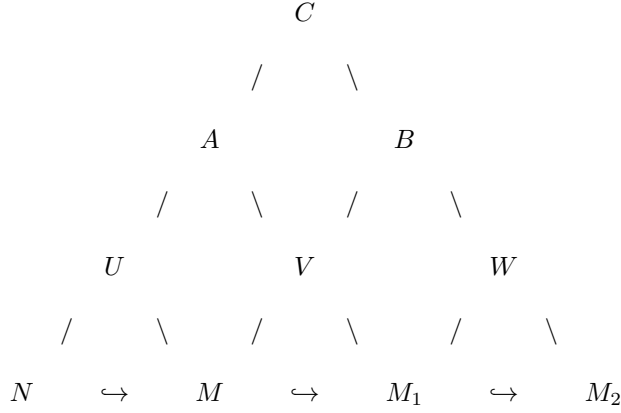

\label{centralizer figure}
$
\begin{array}{rcccccccccccl}
  &     &     &            &   &            & C    &            &     &            &       &           &    \\
  &     &     &            &   &            &     &            &     &            &       &           &    \\
  &     &     &            &   & /          &      &\backslash  &     &            &       &           &    \\
  &     &     &            &   &            &     &            &     &            &       &           &    \\
  &     &     &            & A &            &      &            & B   &            &      &            &   \\
  &     &     &            &   &            &     &            &     &            &       &           &    \\
  &    &      & /          &   &\backslash &       & /          &     & \backslash &      &            &   \\
  &     &     &            &   &            &     &            &     &            &       &           &    \\
  &    & U    &            &   &           &  V    &            &    &             &  W  &           &   \\
  &     &     &            &   &            &     &            &     &            &       &           &    \\
  & / &       & \backslash &   & /          &       & \backslash &    & /           &       & \backslash&   \\
  &     &     &            &   &            &     &            &     &            &       &           &    \\
N &   & \into &            & M &            & \into &            &M_1 &             & \into &           &M_2     
\end{array}$
\caption{Hasse Diagram for Centralizers.}
\end{figure}

\begin{definition}
\label{depth2}
We say that $M/N$ has a (weak) {\em depth $2$} property if the following 
conditions are  satisfied by its Jones tower: 
\begin{enumerate}
\item
$E_M$ has dual bases $\{ z_j \}$, $\{ w_j \} $ in $A$.
\item
$E_{M_1}$ has dual bases $\{ u_i \}$, $\{ v_i \}$ in $B$.
\end{enumerate}
\end{definition}

We note that the depth $2$ conditions in \cite{KN} are a special case
of these.  However, the weak depth $2$ conditions may depend on the choice
of conditional expectation $E: M \rightarrow N$.  

\begin{remark}
If $M/N$ is a subfactor of a finite index von Neumann factor (i.e., $[M:N]<\infty$)
then the above notion of depth $2$ coincides with the usual one. 
\end{remark}

Note that the definition of depth $2$ makes sense for a Frobenius 
extension $M/N$, since for these we retain an endomorphism ring theorem 
stating that Frobenius coordinates $E, x_i, y_i$ for $M/N$ lead to  
coordinates $E_M(m e_1 m') = mm'$ ($m, m' \in M$)
with dual bases 
$\{ x_i e_1 \}$, $\{ e_1 y_i \}$ for $M_1 = M \otimes_N M \cong
\End (M_N)$ as a Frobenius extension over $M$ \cite{O}.  (However,
we no longer necessarily have  $E(1) = 1$ and $e_1^2 = e_1$.)  

We will denote by $T$
the restriction of the normalized trace $T_2 = T_1 E_{M_1}$ of $M_2$ on $C$.

\begin{lemma}
\label{semisimple}
$A$, $B$  are  separable algebras
with $T|_A, T|_B$ as  non-degenerate traces. 
\end{lemma}
\begin{proof}
From the first  of the depth $2$ conditions,
 we see that $E_M(az_j)w_j = a = z_jE_M(w_j a)$ for all
$a \in A \subset M_1$.  Since $z_j w_j = \lambda^{-1} 1$ and $E_M(A) = U$,
we readily see that $A$ is a strongly separable extension of $U$
with Markov trace of index $\lambda^{-1}$.  Similarly,  $B/V$
is a strongly separable extension 
with $E_{M_1}: B \to V$ as conditional expectation, dual bases
$\{ u_i \}$, $\{ v_i \}$ and
index $\lambda^{-1}$.
In particular, $A$ is a separable extension of the separable algebra
$U$,  and is itself a separable algebra \cite{HS}. Similarly,
$B$ is $k$-separable. $T|_A$ is a non-degenerate trace on $A$
since it is a Frobenius homomorphism
by transitivity: $T|_A = T|_U \circ E_M|_A$ by the Markov property.  
Similarly, $T|_B$ is a non-degenerate trace. 
\end{proof}
 
\begin{lemma}
As vector spaces, $M_2 \cong M_1 \otimes_V B$ via the mapping $m_1 \otimes b
\mapsto m_1 b$. Similarly, $M_1 \cong M \otimes_U A$. 
\label{M_2 is M_1B}
\end{lemma}
\begin{proof}
The inverse mapping is given by $x \mapsto E_{M_1}(xu_j) \otimes v_j$.
We note that $$
E_{M_1}(ybu_j) \otimes v_j = y \otimes E(bu_j)v_j = y \otimes b
$$
for $y \in M_1, b \in B$, since $E_{M_1}(B) = V$.   The second
statement can be proven similarly.  
\end{proof}

We develop the following depth $2$ properties for the algebra extension 
$M/N$ above in a series of propositions. 
We let  $E_A = E_{M_1}|_C$. 

\begin{proposition}[Existence of $E_B$]
There exists a $B$-bimodule map $E_B: C\to B$ such that
$E_B|_B = \id_B$,  $E_B$ is a conditional expectation
such that $E_B(e_1) =\lambda 1$
and $T(E_B(c)b) =T(bc)$ for all $b\in B$ and $c\in C$.
\label{Existence of $E_B$}
\end{proposition}
\begin{proof}
Let $\{ a_i \}$, $\{ b_i \}$ denote dual bases in $U$ for $T$
restricted thereon. It follows from Lemma (\ref{V too}) that
 the elements $\{ c_i := \phi(a_i) \}$,
$\{ d_i := \phi(b_i) \}$ are dual bases for the trace
$T$ restricted to $V$.  Define $E_B$ by
\begin{equation}
E_B(c) = T(c u_j c_i)d_iv_j.
\end{equation}
Since $\{u_j c_i \}$, $\{ d_i v_j \}$ are dual bases for $T = TE_{M_1}: B \to k$
by transitivity,
it follows that $E_B(b)= b$ and $E_B(cb) = E_B(c)b$ for every $b \in B$.
The left $B$-module property of $E_B$ follows from: for all $b \in B, c \in C$, 
\[
E_B(bc) = T(bcu_j c_i) d_i v_j = T(cu_j c_i b) d_i v_j = T(cu_jc_i)bd_i v_j,
\]
since $u_j c_i b  \otimes d_i v_j =  u_j c_i  \otimes bd_i v_j$ by Lemma
(\ref{casimir}).

Next,
\[
T(E_B(c)) = T(cu_j c_i) T(d_i v_j) = T(c) 
\]
since $u_j c_iT(d_i v_j) = 1$. 

Finally, let $F = E_M E_A$  and use the Pimsner-Popa relations
as well as the expression for $\phi^{-1}$  below Eq.\ (\ref{phi})
to compute: 
\begin{eqnarray*}
E_B(e_1) = T(e_1u_j c_i)d_i v_j & = & T(e_1 E_A(u_j) c_i) d_i v_j \\
& = & \lambda^{-1} T(e_1 F(e_1u_j) c_i) d_i v_j \\
& = & T(\lambda^{-1} E_M(e_1 c_i) F(e_1 u_j))d_i v_j \\
& = & x_k e_1 T(E_M(e_1(E_A(u_j))a_i)b_i y_k v_j \\
& = & \lambda x_k e_1 E_A(u_j) v_j y_k = \lambda 1_{M_1}. \qed
\end{eqnarray*} 
\renewcommand{\qed}{}\end{proof} 

\begin{proposition}[``Commuting square condition'']
\label{commuting square condition}
 We have $E_A\circ E_B = E_B\circ E_A$.
\end{proposition}
\begin{proof}
We compute: for each $c \in C$, 
\[
E_A E_B(c) = T(cu_j c_i) d_i E_A(v_j) = T(cu_j E_A(v_j)c_i)d_i= T(cc_i)d_i,
\]
while
\[
E_BE_A(c) =  T(E_A(c)E_A(u_j)c_i) d_i v_j = 
T(E_A(c) c_i) d_i E_A(u_j) v_j = T(cc_i)d_i
\] 
by the Markov property $T E_A = T$. 
\end{proof}

\begin{proposition}[``Symmetric square condition'']
\label{symmetric square condition}
We have $AB = BA =C$. More precisely, 
$A\otimes _V B \cong B\otimes_V A \cong C$ as vector spaces. 
\end{proposition}
\begin{proof}
We note that $E_A(C) = A$ and $V = A \cap B$. 
The proposition follows easily from the dual bases equations and the depth $2$ assumption:  
\[
 E_A(cu_j) v_j = c = u_j E_A(v_j c),
\]
for all $c \in C$.  
\end{proof} 

\begin{proposition}[Pimsner-Popa identities]
We have
\begin{eqnarray*}
\lambda^{-1} e_2 E_A(e_2c) &=& e_2c, \qquad
\lambda^{-1} E_A(ce_2)e_2  = c e_2, \\
\lambda^{-1} e_1 E_B(e_1c) &=& e_1c, \qquad
\lambda^{-1} E_B(ce_1)e_1  = c e_1.
\end{eqnarray*}
As a consequence we have 
\begin{eqnarray*}
Ce_2 &=& Ae_2, \qquad e_2C = e_2A, \\
Ce_1 &=& Be_1, \qquad e_1C = e_1B. 
\end{eqnarray*}

\end{proposition}
\begin{proof}
Now $e_2C = e_2 A$ and $Ce_2 = Ae_2$ follow from the usual Pimsner-Popa Eqs.\
(\ref{PP}) 
for $E_{M_1}|_C = E_A$. At a point below in this proof, we will need to know that 
\begin{equation}
C = Ae_2 A,
\label{C basic over A}
\end{equation}
which follows from
\[
c = E_A(cu_j) v_j = \lambda^{-1} E_A(cz_ie_2) e_2 w_i,
\]
for by the basic construction theorem 
$u_j \otimes v_j = \lambda^{-1} z_i e_2 \otimes e_2 w_i$
in $M_2 \otimes_{M_1} M_2$. 

Note that $F(C) = U$. 
We compute:  for each $c \in C$,
\begin{eqnarray*}
e_1 c = e_1 E_A(cu_j) v_j & = & \lambda^{-1} e_1 E_M(e_1 E_A(cu_j))v_j \\
& = & \lambda^{-1} e_1 T(F(e_1 c u_j)a_i)b_i v_j \\
& = & \lambda^{-3} T(cu_j E_M(c_ie_1)e_1) e_1 E_M(e_1 d_i) v_j \\
& = & \lambda^{-1} T(cu_j c_i e_1) e_1 d_i v_j = \lambda^{-1} e_1 E_B(e_1c).
\end{eqnarray*}
Thus, $e_1 C = e_1 B$. 

The computation $ce_1 = \lambda^{-1} E'_B(ce_1)e_1$ proceeds similarly, 
 where
\begin{equation}
E'_B(c) = u_j c_i T(d_i v_j c),
\end{equation}
clearly defines a bimodule projection of $C$ onto $B$ 
(cf.\ Proposition (\ref{Existence of $E_B$})).  
As a result, we have $Ce_1 = Be_1$. 

We will show that $E_B = E'_B$ by showing that $C = Be_1B$ and 
noting that $E'_B(e_1) = \lambda 1$ by a computation very similar to that for 
$E_B(e_1) = \lambda 1$ above. Using the braid-like relations and Eq.\ (\ref{C basic over A}),
we compute:
\[
C = Ae_2 A = Ae_2 e_1 e_2 A = Ce_1C = Be_1 B. \qed \] 
\renewcommand{\qed}{}\end{proof}
 
It is not hard to show that $E_B: C \to B$  is isomorphic to the
basic construction of the strongly separable extension $B/V$, where $C = Be_1 B$.
Similarly, $E_A: C \to A$ is isomorphic to the basic construction of
the strongly separable extension $A/U$, where $C = Ae_2A$. 

\begin{remark}
Irreducible separable Markov extensions considered in \cite{KN}
trivially satisfy the weak irreducibility assumption as well as the conclusion
of Lemma (\ref{V too}). 
It follows that all the results of the next sections apply to these. 
\end{remark}
 
\end{section}


\begin{section}
{Weak Hopf algebra structures on centralizers}

 Let $f =f^{(1)}\otimes f^{(2)}$
be the unique symmetric separability element \cite{SK} of $V = C_{M_1}(M)$,
and let $w = [f^{(1)} T(f^{(2)})]^{-1} \in Z(V)$ be the invertible element
satisfying $f^{(1)} T(vwf^{(2)}) = v$ for all $v\in V$. In other words,
$f^{(1)}\otimes wf^{(2)}$ is the dual bases tensor for  $T: V \to k$.

\begin{proposition}
\label{duality}
The  following bilinear form,
$$
\la a, \, b \ra = \lambda^{-2} T(ae_2e_1wb), \qquad a\in A,\, b\in B,
$$
is non-degenerate on $A\otimes B$.
\end{proposition}
\begin{proof}
If $ \la a,\, B\ra =0$ for some $a\in A$, then
for all $x\in C$ we have $T(ae_2e_1x)=0$,
since $e_1B = e_1 C$ (depth $2$ property).
Taking $x = e_2a'\,(a'\in A)$ and using 
the braid-like relation between
Jones idempotents, and Markov property of $T$ we have 
$$
T(aa') = \lambda^{-1}T(ae_2e_1(e_2a')) = 0\qquad \mbox{ for all } a'\in A,
$$ 
therefore $a=0$. Similarly, one proves that $\la A,\, b\ra =0$
implies $b=0$. 
\end{proof}

The above duality form allows us to introduce a comultiplication 
$b\mapsto b\1 \otimes b\2$ on $B$ as follows:
\begin{equation}
\la a_1,\, b\1 \ra \la a_2,\, b\2 \ra = \la a_1a_2,\, b\ra,
\label{comultiplication}
\end{equation} 
for all $a_1,a_2\in A,\, b\in B$, and counit $\eps: B \rightarrow k$
by ($\forall b \in B$)
\begin{equation}
\eps(b) = \la 1,\,b \ra.
\label{counit}
\end{equation} 

A proof similar to that of Proposition (\ref{duality})   shows that 
$\la a,b \ra' = \lambda^{-2} T(be_1e_2 wa)$ is
another non-degenerate pairing of $A$ and $B$. 
We then introduce a linear automorphism $S: B \to B$ by the relation
$\la a,b \ra = \la a,S(b) \ra'$, i.e.,
\begin{equation}
\la a, b \ra  =\lambda^{-2} T(S(b)e_1e_2wa)
\end{equation} 
for all $a\in A,\, b\in B$, or, equivalently,
\begin{equation}
\label{eqn : S}
E_A(e_2e_1wb) = E_A(S(b) e_1e_2)w. 
\end{equation}
Note that we automatically have
\begin{equation}
E_{M_1}(e_2xwb) = E_{M_1}(S(b)xe_2)w, \qquad \mbox{for all } x\in M_1.
\end{equation} 

\begin{proposition}
\label{1 and eps}
We note that: (for all $b,c\in B$)
\begin{equation}
\eps(b) = \lambda^{-1}T(e_2wb),
\end{equation}
\begin{equation}
\eps(S(b)) = \eps(b),
\end{equation}
\begin{equation}
\Delta(1) = S^{-1}(f^{(1)})\otimes f^{(2)}.
\end{equation}
\end{proposition}
\begin{proof}
The formula for $\eps$ follows from the identity $E_B(e_1) = \lambda 1$ 
and $T\circ E_B = T$:
\begin{equation*}
\eps(b) = \lambda^{-2} T(e_2e_1wb) = \lambda^{-1}T(e_2wb).
\end{equation*}
Then the second equation follows from the computation :
$$
\eps(b) = \lambda^{-1} T(e_2 wb) = \lambda^{-2} T(b E_B(e_1) e_2 w)
= \lambda^{-2} T(e_2 e_1 w S^{-1}(b)) = \eps(S^{-1}(b)).  
$$
To establish the third formula, we use the
Markov property and commuting square condition to compute: for all $a,a'\in A$, 
\begin{eqnarray*}
\la a,  S^{-1}(f^{(1)})\ra \la a',\,f^{(2)} \ra
&=&  \lambda^{-3} T(ae_2e_1w S^{-1}(f^{(1)})) T(E_A\circ E_B(a'e_1w) f^{(2)})\\
&=& \lambda^{-3} T(f^{(1)}e_1e_2wa) T(E_B(a'e_1) w f^{(2)}) \\
&=& \lambda^{-2} T(E_B(a'e_1)e_1e_2wa) \\
&=& \lambda^{-2} T(aa'e_1e_2w) =\la aa',\,1\ra. \qed 
\end{eqnarray*}
\renewcommand{\qed}{}\end{proof} 

The following lemma gives a useful explicit formula for $S^{-1}$.

\begin{lemma}
\label{explicit S}
For all $b\in B$ we have  
$S^{-1}(b)= \lambda^{-3} w^{-1} E_B(e_1e_2 E_A(be_1e_2)) w$.
\end{lemma}
\begin{proof}
We obtain this formula by multiplying  both sides of
Eq.\ (\ref{eqn : S}) by $e_1e_2$ on the left 
and taking $E_B$  from both sides.
\end{proof}

\begin{corollary}
\label{S maps bases}
We have $S(V) = W$, where $W= C_{M_2}(M_1)$.
\end{corollary}
\begin{proof}
Let us take $y\in W$, then using Lemma (\ref{explicit S}), the commuting square
condition, 
and the Markov property we have 
\begin{eqnarray*}
S^{-1}(y)
&=& \lambda^{-3} w^{-1} E_B(e_1e_2 e_1 E_A(ye_2)) w \\
&=& \lambda^{-2} w^{-1} E_B(e_1 E_A(ye_2)) w \in V.
\end{eqnarray*}
Therefore, $S^{-1}(W) \subseteq V$ and since $W\cong V$ as vector spaces,
we have  $S(V) = W$.
\end{proof}

\begin{lemma}
\label{toward S2}
For all $b\in B$ we have  
$b = w S^{-1}(wS^{-1}(b)w^{-1}) w^{-1}$.
\end{lemma}
\begin{proof}
Using non-degeneracy of the duality form and definition of $S$ we compute
for all $a\in A$:
\begin{eqnarray*}
T(ae_2e_1b)
&=& \lambda^{-1} T(E_A(bae_2)e_2e_1) \\
&=& \lambda^{-1} T(E_A(e_2awS^{-1}(b))w^{-1}e_2e_1) \\
&=& T(e_2awS^{-1}(b)w^{-1} e_1) \\
&=& T(E_A(wS^{-1}(b)w^{-1}e_1e_2) w w^{-1} a) \\
&=& T(a E_A(e_2e_1wS^{-1}(wS^{-1}(b)w^{-1}))w^{-1}),
\end{eqnarray*}
whence the formula follows.
\end{proof}

\begin{proposition}
\label{anti algebra}
$S$ is an algebra anti-homomorphism, i.e.,
\begin{equation*}
S(bb') =S(b')S(b) \qquad \mbox{for all } b,b'\in B.
\end{equation*}
\end{proposition}
\begin{proof}
We use the non-degeneracy of the duality form:
\begin{eqnarray*}
T(ae_2e_1wS^{-1}(b')w^{-1}S^{-1}(b))
&=& \lambda^{-1} T(w^{-1}E_A(S^{-1}(b)ae_2)e_2e_1wS^{-1}(b')) \\
&=& \lambda^{-1} T(E_A(w^{-1}e_2awS^{-2}(b))w^{-1}e_2e_1wS^{-1}(b')) \\
&=& \lambda^{-1} T(b'e_1e_2 E_A(e_2awS^{-2}(b))w^{-1}) \\
&=& T(wS^{-2}(b)w^{-1}b'e_1e_2aw) \\
&=& T(ae_2e_1w S^{-1}(wS^{-2}(b)w^{-1}b')w^{-1}),
\end{eqnarray*}
therefore, we have $S^{-1}(b')w^{-1}S^{-1}(b)w =  S^{-1}(wS^{-2}(b)w^{-1}b')$.
Using Lemma (\ref{toward S2}) we conclude that 
\begin{equation*}
S^{-1}(b') S^{-1}(w S^{-2}(b) w^{-1}) = S^{-1}(b') w^{-1} S^{-1}(b) w 
= S^{-1}(w S^{-2}(b) w^{-1} b').
\end{equation*}
We replace $wS^{-2}(b)w^{-1}$ by $b$ to obtain 
the result.
\end{proof}

\begin{corollary}
\label{S2 on bases}
For all $b\in B$ we have  $S^2(b) = gbg^{-1}$ where 
\begin{equation}
\label{g}
g := S(w^{-1})w.
\end{equation}
\end{corollary}

In particular, $S^2|_V = \id_V$ from (\ref{S maps bases}), 
so $S$ maps $V$ to $W$ and vice versa, as well
as $S^2|_{W}= \id_{W}$.
For example, we obtain $\Delta(1) = S(f^{(1)}) \otimes f^{(2)}$ from this
and (\ref{1 and eps}). 

\begin{lemma}
\label{delta on V}
For all $b\in B$ and $v\in V$ we have
\begin{equation}
\Delta(bv) = \Delta(b) (v\otimes 1).
\end{equation}
\end{lemma}
\begin{proof}
Let $a,a'\in A$ then 
\begin{eqnarray*}
\la a\otimes a',\, \Delta(bv) \ra 
&=& \la aa',\, bv\ra = \la vaa',\, b\ra \\
&=& \la a,\,b\1v\ra \la a',\,b\2 \ra. \qed 
\end{eqnarray*}
\renewcommand{\qed}{}\end{proof}

Now we are in the position to establish the unit and counit
axioms for $B$.

\begin{proposition}
\label{counit axioms}
We have 
\begin{equation}(\id\otimes\Delta)\Delta(1) = 
(\Delta(1)\otimes 1)(1\otimes \Delta(1)) = 
(1\otimes \Delta(1))(\Delta(1)\otimes 1).
\end{equation}
\end{proposition}
\begin{proof}
We have seen that $\Delta(1)\in W\otimes V$, therefore
$(1\otimes \Delta(1))$ and $(\Delta(1)\otimes 1)$ commute.
By Lemma (\ref{delta on V}),
\begin{eqnarray*}
(1\otimes \Delta(1))(\Delta(1)\otimes 1)
&=& S^{-1}(f^{(1)}) \otimes 1\1 f^{(2)} \otimes 1\2 \\
&=& S^{-1}(f^{(1)}) \otimes \Delta(f^{(2)}) 
    = (\id\otimes\Delta)\Delta(1). \qed 
\end{eqnarray*}
\renewcommand{\qed}{}\end{proof}

\begin{proposition}
\label{axiom for counit}
For all $b,c,d\in B$ we have 
\begin{equation*}
\eps(bcd) =\eps(bc\1)\eps(c\2d) =  \eps(bc\2)\eps(c\1d).
\end{equation*}
\end{proposition}
\begin{proof}
First, one can define a coalgebra structure on $A$ using the
duality form from Proposition (\ref{duality}) and show that
$\Delta_A(1)\in A \otimes C_M(N)$ in a way similar to how it
was shown above for the comultiplication $\Delta$ of $B$ 
that $\Delta(1) \in W \otimes V$. Then  we compute:
\begin{eqnarray*}
\eps(bcd)
&=& \lambda^{-1} T(e_2wbcd) \\
&=& \lambda^{-3} T(E_A(de_2)e_2e_1 wbc) \\
&=& \la 1\1,\, b \ra  \la \lambda^{-1} E_A(de_2)1\2,\,c \ra \\
&=& \la 1\1,\, b \ra \la 1\2,\, c\2 \ra  \la \lambda^{-1} E_A(de_2),\,c\1 \ra\\
&=& \eps(bc\2) \eps(c\1d).
\end{eqnarray*}
Note that in the third line $E_A(de_2)$ commutes with each of the
elements in $\{ 1\2 \} \subset U$, 
so that $\eps(bcd)$ is also equal to $\eps(bc\1)\eps(c\2d)$.
\end{proof}

The next step is to prove that $\Delta$ is a homomorphism.
To achieve this we first need to establish a certain commutation relation
(see Proposition (\ref{heart}) below) that corresponds to the two
different ways of representing $C=AB =BA$.

We will need several preliminary results.

\begin{lemma}
\label{abcde}
The following identities hold for all $b\in B$ and $v\in V$:
\begin{enumerate}
\item[(a)]  $S^{-1}(e_2) = w^{-1}e_2w$,
\item[(b)] $ve_2 = S(v)e_2$,
\item[(c)] $\lambda^{-1} E_A(e_2wb)w^{-1} = \eps(b1\1)1\2$,
\item[(d)] $\Delta(b) (1\otimes v) = \Delta(b) (S(v)\otimes 1)$,
\item[(e)] $\Delta(b)\Delta(1) = \Delta(b)$.
\end{enumerate}
\end{lemma}
\begin{proof}
(a) We have $T(ae_2e_1wS^{-1}(e_2) ) = T(e_2e_1e_2wa) = T(ae_2e_1e_2w)$,
whence the result follows by non-degeneracy of the bilinear pairing 
$a\otimes b \mapsto T(ae_2e_1b)$.

(b) We compute, using part (a) and the anti-multiplicativity of $S$: 
\begin{eqnarray*}
\lambda^2 \la a,\, S^{-1}(ve_2) \ra
&=& T(ve_2e_1e_2wa) \\
&=& T(ae_2e_1wS^{-1}(ve_2)) \\
&=& \lambda T(ae_2wS^{-1}(v)) \\
&=& T(S^{-1}(v)e_2e_1e_2wa) = \lambda^2 \la a,\, S^{-1}(S^{-1}(v)e_2) \ra.
\end{eqnarray*}

(c) Since both sides of the given equation belong to $V$,  it suffices
to evaluate them against $T(\cdot\, v)$ for all $v\in V$:
\begin{eqnarray*}
T(\lambda^{-1} E_A(e_2wb)v) 
&=&  \lambda^{-1} T(e_2wbv) = \lambda^{-1}T(ve_2wb) \\
T(\eps(b1\1)1\2wv)
&=& \eps(bS(f^{(1)}))T(vwf^{(2)}) \\
&=& \eps(bS(v)) = \lambda^{-1} T(e_2wbS(v)) \\
&=& \lambda^{-1}T(ve_2wb),
\end{eqnarray*}
where we used part (b).

(d) We evaluate both sides against elements $a\otimes a' \in A\otimes A$ 
(note that $S(v)$ commutes with $A$):
\begin{eqnarray*}
\la a\otimes a',\, b\1S(v)\otimes b\2 \ra
&=& \lambda^{-2} T(S(v)ae_2e_1wb\1)\la a',\,b\2 \ra \\
&=& \lambda^{-2} T(ave_2e_1wb\1) \la a',\,b\2 \ra \\
&=& \la av,\, b\1 \ra \la a',\, b\2\ra = \la ava',\, b\ra \\
&=& \la a\otimes a',\, b\1 \otimes b\2v  \ra.
\end{eqnarray*}

(e) From part (d), properties of $S$, 
and properties of the separability element $f$ we have
\begin{eqnarray*}
\Delta(b)\Delta(1) 
&=& b\1 1\1 \otimes b\2 1\2 = b\1 S(1\2)1\1 \otimes b\2 \\
&=& b\1 S(f^{(1)} f^{(2)}) \otimes b\2 = b\1 \otimes b\2. \qed 
\end{eqnarray*}
\renewcommand{\qed}{}\end{proof}

Applying $S$ to part (a) above, we obtain from part (b):
\begin{equation}
S(e_2) = w^{-1} e_2 w.
\end{equation}

\begin{proposition}
\label{two identities}
For all $a\in A$ and $b\in B$ we have
\begin{enumerate}
\item[(i)] $ \lambda^{-1}E_B(e_1wba) = \la a,\, b\1 \ra w b\2$, 
\item[(ii)] $ \lambda^{-1}b\2 E_A(e_2wb\1) w^{-1} = b$.
\end{enumerate}
\end{proposition}
\begin{proof}
(i) Let $a'\in A$ then
\begin{eqnarray*}
\la a',\, \lambda^{-1} w^{-1}E_B(e_1wba) \ra
&=& \lambda^{-3} T(a'e_2e_1 E_B(e_1wba')) \\
&=& \lambda^{-2} T(a'e_2e_1wba') = \la aa',\, b\ra \\
&=& \la a',\, \la a,\, b\1 \ra b\2 \ra.
\end{eqnarray*}

(ii) From Lemma (\ref{abcde}c) and (e) we have
\begin{equation*}
\lambda^{-1}b\2 E_A(e_2wb\1) w^{-1} = \eps(b\1 1\1) b\2 1\2 = b. \qed
\end{equation*}
\renewcommand{\qed}{}\end{proof}

The next Proposition (cf.\ \cite{KN}, 4.6) is the key 
ingredient in proving that $B$ is a weak Hopf algebra
acting on $M_1$.

\begin{proposition}
\label{heart}
For all $b\in B$ we have 
\begin{equation}
w^{-1}e_1wb = \lambda^{-1} b\2 w^{-1} E_A(e_2e_1wb\1).
\end{equation}
\end{proposition} 
\begin{proof}
First, let us note that for all $c_1, c_2 \in C$ we have
$c_1 =c_2$ if and only if $E_B(c_1a) =E_B(c_2a)$ for all $a\in A$.
Indeed, if $c\in C$ and $E_B(ca)=0$ for all $a\in A$ then 
$T(abc) = T(b E_B(ca)) =0$ for all $b\in B$. But since $AB =C$
by Proposition (\ref{symmetric square condition})
and $T$ is non-degenerate, we conclude that $c=0$.

Let $c_1 = w^{-1}e_1wb$ and $c_2 = \lambda^{-1} b\2 w^{-1} E_A(e_2e_1wb\1)$.
We compute, using Propositions (\ref{two identities}) and 
(\ref{commuting square condition}) : 
\begin{eqnarray*}
E_B(c_1a)
&=& w^{-1}E_B(e_1wba) = w^{-1} \la a,\, b\1 \ra w b\2 \\
&=& \la a,\, b\1 \ra b\2, \\
E_B(c_2a)
&=& \lambda^{-1} b\2 w^{-1} E_B\circ E_A(e_2e_1wb\1 a) \\
&=& \lambda^{-1} b\2 w^{-1} E_A(e_2 E_B(e_1wb\1 a)) \\
&=& \lambda^{-1} \la a,\, b\1 \ra b\3 w^{-1}  E_A(e_2wb\2)\\
&=& \la a,\, b\1 \ra b\2,
\end{eqnarray*}
whence the result follows.
\end{proof}

\begin{corollary} 
\label{M1 corollary}
For all $b\in B$ and $x\in M_1$ we have
\begin{equation}
w^{-1}xb = \lambda^{-1} b\2 w^{-1} E_{M_1}(e_2xb\1).
\label{full equation}
\end{equation}
\end{corollary}
\begin{proof}
This follows from the fact that every $x\in M_1$ can be written
as $x= \sum\, x_ie_1y_i$, where $x_i, y_i\in M$ commute with $B$.
\end{proof}

\begin{corollary} 
\label{action corollary}
For all $x,y\in M_1$ and $b\in B$, we have
\begin{equation}
E_{M_1}(e_2wyxb) =  \lambda^{-1} E_{M_1}(e_2wyb\2) w^{-1} E_{M_1}(e_2wxb\1).
\label{anti-measuring}
\end{equation}
\end{corollary}
\begin{proof}
This is obtained from Corollary (\ref{M1 corollary}) by replacing
$x$ with $wx$, multiplying both sides by $e_2wy$ on the left,
and taking $E_A$ from both sides.
\end{proof}

In order to prove the multiplicativity of $\Delta$ we first need
to establish anti-comultiplicativity of $S$.

\begin{proposition}
\label{anti coalgebra}
$S$ is anti-comultiplicative, i.e., 
\begin{equation}
\Delta S(b) = S(b\2) \otimes S(b\1) \qquad \mbox{for all } b\in B.
\end{equation}
\end{proposition}
\begin{proof}
Let $a,a'\in A$ then using Corollary (\ref{action corollary})
and Lemma (\ref{abcde}d) we compute:
\begin{eqnarray*}
\la aa',\, S^{-1}(b) \ra
&=& \lambda^{-3} T(e_1e_2 E_A(e_2waa'b)) \\
&=& \lambda^{-4} T(e_1e_2 E_A(e_2wab\2)w^{-1}E_A(e_2wa'b\1)) \\
&=& \lambda^{-2} \la w^{-1} E_A(e_2wab\2)w^{-1} E_A(e_2wa'b\1),\, 1\ra \\
&=& \lambda^{-2} \la w^{-1} E_A(e_2wab\2),\, 1\1 \ra
                 \la w^{-1} E_A(e_2wa'b\1),\, 1\2 \ra \\
&=& \lambda^{-6} T(S(1\1)e_1e_2E_A(e_2wab\2)) T(S(1\2)e_1e_2E_A(e_2wa'b\1))\\
&=& \lambda^{-4} T(b\2 S(1\1) e_1e_2wa) T(b\1 S(1\2) e_1e_2wa') \\
&=& \la a,\, S^{-1}(b\2S(1\1))\ra  \la a',\, S^{-1}(b\1S(1\2))\ra \\
&=& \la a,\, S^{-1}(b\2)\ra \la a',\, S^{-1}(b\1 1\1 S(1\2)) \ra \\
&=& \la a,\, S^{-1}(b\2)\ra \la a',\, S^{-1}(b\1)\ra,
\end{eqnarray*}
since $f^{(2)}f^{(1)} = 1$,
whence the proposition follows from non-degeneracy of $\la\, ,\, \ra$
and bijectivity of $S$.
\end{proof}

\begin{proposition}
\label{multiplicativity}
$\Delta$ is a homomorphism of algebras:
\begin{equation}
\Delta(bb') = \Delta(b) \Delta(b') \qquad \mbox{for all } b,b'\in B.
\end{equation}
\end{proposition}
\begin{proof}
Using the definition and properties of $S$ and 
Corollary (\ref{action corollary}) for all $x,y \in M_1$ we have:
\begin{eqnarray*}
E_{M_1}(S(b)xw^{-1}ye_2)w
&=& E_{M_1}(e_2xw^{-1}ywb) \\
&=& \lambda^{-1} E_{M_1}(e_2xb\2) w^{-1} E_{M_1}(e_2ywb\1) \\
&=& \lambda^{-1} E_{M_1}(S(b\2)xw^{-1}e_2) E_{M_1}(S(b\1)ye_2)w,
\end{eqnarray*}
and using Corollary (\ref{anti coalgebra}) and bijectivity of $S$
we obtain:
\begin{equation}
\label{baa'}
E_{M_1}(bxye_2) = \lambda^{-1} E_{M_1}(b\1 x e_2) E_{M_1}(b\2 y e_2)
\qquad \mbox{for all } x,y \in M_1, b\in B.
\end{equation}
Next, using the duality form we have: for  $a,a' \in A$, 
\begin{eqnarray*}
\la a\otimes a',\, \Delta(bb') \ra
&=& \la aa',\, bb'\ra \\
&=& \lambda^{-1} \la E_A(b'aa'e_2), b \ra \\
&=& \lambda^{-2} \la E_A(b\1'ae_2), b\1 \ra \la  E_A(b\2'ae_2), b\2 \ra\\
&=& \la a,\, b\1b\1'\ra   \la a',\, b\2b\2'\ra, 
\end{eqnarray*}
as required.
\end{proof}

Next we establish properties of the antipode with respect to
the counital maps.

\begin{proposition}
\label{prop:counital maps}
For all $b\in B$ we have the following identities:
\begin{eqnarray}
S(b\1)b\2 &=& 1\1 \eps(b 1\2), \\
b\1 S(b\2) &=& \eps(1\1b)1\2.  \label{left counital formula}
\end{eqnarray}
\end{proposition}
\begin{proof}
To establish the first relation we compute, using Eq.\ (\ref{baa'}),
for all $a\in A$:
\begin{eqnarray*}
\la a,\, S^{-1}(b\1)w^{-1}b\2 \ra
&=& \lambda^{-1} \la E_A(w^{-1}b\2ae_2),\, S^{-1}(b\1) \ra \\
&=& \lambda^{-4} T( E_A(w^{-1}b\2ae_2) e_2 E_A(e_2e_1w S^{-1}(b\1))) \\
&=& \lambda^{-3} T( E_A(b\2ae_2) E_A(b\1e_1e_2) ) \\
&=& \lambda^{-2} T(be_1a e_2).
\end{eqnarray*}
Next we recall the formula for $\Delta(1)$ from Proposition (\ref{1 and eps}),
formula for $S^2$ from Corollary (\ref{S2 on bases}), 
Lemma (\ref{abcde}d), and that $\Delta(w) =\Delta(1)(w\otimes 1) 
= (w\otimes 1)\Delta(1)$:
\begin{eqnarray*}
\la a,\, 1\1\eps(b1\2) \ra
&=& \lambda^{-1} \la a,\, 1\1\ra T(e_2wb1\2) \\
&=& \lambda^{-1} \la a,\, S^{-1}(f^{(1)})\ra T(e_2wb f^{(2)}) \\
&=& \lambda^{-1} \la a,\, S^{-1}(E_A(e_2wbw^{-1})) \ra \\
&=& \lambda^{-3} T( E_A(e_2wbw^{-1})e_1e_2wa) \\
&=& \lambda^{-2} T(e_2wbw^{-1}e_1wa) \\
&=& \la wa,\, S^{-1}((wbw^{-1})\1)  w^{-1} (wbw^{-1})\2 \ra \\
&=& \la a,\, S^{-1}(wb\1w^{-1}) w^{-1} b\2 w \ra \\
&=& \la a,\, S^{-1}(w S(w^{-1}) b\1 S(w) w^{-1} ) b\2 \ra \\
&=& \la a,\, S(b\1)b\2 \ra.
\end{eqnarray*}
The second identity follows from the first by (\ref{1 and eps}),
since the symmetry of $f$ and the
anti-(co)multiplicative properties of the antipode imply :
\begin{eqnarray*}
b\1 S(b\2)
&=& S( S(S^{-1}(b)\1) S^{-1}(b)\2) = S(1\1) \eps (S^{-1}(b) 1\2) \\
&=& \eps (S(1\2)b) S(1\1) = \eps(1\1b)1\2. \qed 
\end{eqnarray*}
\renewcommand{\qed}{}\end{proof}

Let us consider two mappings $\eps_t: B \to V$
and $\eps_s: B \to W$ given by $\eps_t(b) = \eps(1\1b)1\2 $,
 and $\eps_s(b) = 1\1 \eps(b 1\2)$, 
corresponding to the right-hand side of the equations in Proposition (\ref{prop:counital maps}).  They are  called the \textit{target and source counital maps}, 
respectively (cf.\ Section~1). 
By a computation quite similar to that in Lemma (\ref{abcde}c),
we may check that:
\begin{equation}
\eps_t(b) = \lambda^{-1} E_A(be_2).
\label{eps^L}
\end{equation}
Indeed, we have for each $v \in V$, 
$$ T(\eps(1\1 b) 1\2 v) = \eps( S(vw^{-1}) b) = \lambda^{-1} T(e_2 w S(w^{-1}) S(v) b) = \lambda^{-1} T(e_2 vb)
$$
while also $T(\lambda^{-1} E_A(be_2) v) = \lambda^{-1} T(e_2 vb)$.

\begin{theorem}
$(B, \Delta, \eps, S)$ is a semisimple weak Hopf algebra.
\end{theorem}
\begin{proof}
Semisimplicity follows from Lemma (\ref{semisimple}). 
We have established all the axioms of a weak Hopf algebra
 except Axiom (\ref{S id S}), which we
show next.
At a point below, we let $b' = S(b)$, at another $b'' = wb'$, and  use 
Eq.\ (\ref{full equation}) as well as 
Lemma (\ref{delta on V}). Let $g = S(w^{-1})w$ 
be the element from Corollary (\ref{S2 on bases})
implementing the inner automorphism $S^2$, then
for all $b \in B$, 
\begin{eqnarray*}
S(b\1) b\2 S(b\3) & = &  \lambda^{-1} S(b\1) E_A(b\2 e_2) \\
& = & \lambda^{-1} b'\2 E_A(S^{-1}(b'\1)e_2) \\
& = & \lambda^{-1} b'\2 E_A(e_2 w g^{-1} b'\1 g) w^{-1} \\
& = & \lambda^{-1} b'\2 E_A(e_2 w b'\1 S(w^{-1})) \\
& = & \lambda^{-1} b''\2 w^{-1} E_A(e_2 b''\1) = w^{-1} b'' = S(b). \qed
\end{eqnarray*}
\renewcommand{\qed}{}\end{proof}



\begin{remark}
\begin{enumerate}
\item[(i)] $V= \eps_t(B)$ is the target counital subalgebra of $B$
and $W = C_{M_2}(M_1) = S(V)$ is the source counital subalgebra
(recall that the antipode maps one counital subalgebra to another).
\item[(ii)]
From Eq.\ (\ref{eps^L}) we see that $e_2$ is
a normalized left integral in $B$:
$$
be_2 = \lambda^{-1} E_A(be_2)e_2 = \eps_t(b) e_2.
$$
Furthermore, $l = e_2 S^{-1}(e_2) = e_2w^{-1}e_2w = e_2w$
is a two-sided integral in $B$, due to Lemma (\ref{abcde}a)
and the fact that the space of left (respectively, right)
integrals in a weak Hopf algebra is a left (respectively, right)
ideal. Next, $S(l) = w^{-1} S(w) e_2 w = e_2 w =l$, since
$\eps_t|_W = S|_W$.  Finally $l$ is normalized, since
\begin{equation*}
\eps_t(l) = \lambda^{-1} E_A(E_M(w)e_1)=1 \quad
\mbox{ and } \quad \eps_s(l)= S\circ\eps_t(l) =1.
\end{equation*}
Clearly, $l$ is the unique element with these properties
(cf.\ \cite[5.7]{NV1}). Such a two-sided normalized integral
is called a {\em Haar integral} in \cite{BNSz}.
\end{enumerate}
\end{remark}

Defining a comultiplication  and counit  of $A$ similarly to 
Eqs.\ (\ref{comultiplication}) and (\ref{counit}), as the dual
of the multiplication and unit of $B$, and an antipode $S_A$ on $A$ by
$
\la S_A(a), b \ra = \la a, S(b) \ra
$, 
the corollary below
follows from the self-duality of the axioms of a weak Hopf algebra
and Lemma (\ref{semisimple}). 
\begin{corollary}
$A$ is a semisimple weak Hopf algebra isomorphic to the  dual of $B$. 
\end{corollary}  

\end{section}
 \begin{section}{Action and smash product}

In this section we define an action of $B$ on $M_1$
suggested by the measuring in Eq.\ (\ref{baa'}), and show 
that it comes from the standard left action of a weak Hopf 
algebra on its dual. 
 We then show that $M$ is the subalgebra of
invariants of this action, and  that $M_2$ is isomorphic to the 
smash product of $M_1$ with $B$.

\begin{proposition}
\label{action of B on M_1}
 The mapping $\lact:  B \otimes M_1 \to M_1 $ given by
\begin{equation}
b \lact x = \lambda^{-1} E_{M_1}(bxe_2)
\label{act}
\end{equation}
defines a left action of a weak Hopf algebra on $M_1$, characterized by 
\begin{equation}
b \lact ma = m \la a\2, b \ra a\1
\label{standard action}
\end{equation}
for all $m \in M, a \in A, b \in B$.
In particular, $M$ is the subalgebra of invariants for this action.
\end{proposition}
\begin{proof}
From Eq.\ (\ref{baa'}) it follows that $\lact$ satisfies the measuring
axiom.  From Eq.\ (\ref{eps^L}) it follows that
$b \lact 1 = \eps_t(b)$.  The action  of $B$ on $M_1$ is a left module
action of an algebra by the Pimsner-Popa relations and $E_{M_1}(xe_2) = 
\lambda x$ for $x \in M_1$. 

Recall that $M_1 = MA$. Since $B = C_{M_2}(M)$, it is clear that 
$b \lact ma = m (b \lact a)$ for every $m \in M$.  
We compute for every $a \in A, b,b' \in B$:
$$
\la a\1, b' \ra \la a\2, b \ra  =  \la a, b' b \ra 
      =  \la  \lambda ^{-1} E_A(bae_2),b' \ra 
       =  \la b \lact a, b' \ra,
$$
whence Eq.\ (\ref{standard action}) follows.
Thus the action of $B$ on $A$ coincides
with the standard left action of a weak Hopf algebra
$B$ on its dual $B^* \cong A$ as in Example (\ref{examples of actions}(ii)).
Since the invariant subalgebra $A^B$ is  $k 1$, it follows that
$M_1^B = M$.  
\end{proof}

The next proposition provides a simplifying formula for this
action.  We will need  the equation
\begin{equation}
b\1 S(b\2) b\3 = b
\label{WHA theory}
\end{equation}
for each $b \in B$, which follows from Eq.\ (\ref{left counital formula}). 

\begin{proposition}
\label{simplification}
For every $b \in B, x \in M_1$, we have
$$
b \lact x = b\1 x S(b\2).
$$
\end{proposition}
\begin{proof}
We subsequently use Eq.\ (\ref{full equation}), Lemma (\ref{abcde}d)
and its opposite (obtained by applying $S \otimes S$), 
Corollary (\ref{S2 on bases}), and Eq.\ (\ref{WHA theory})
in the next computation:  for every $b \in B, x \in M_1$,
\begin{eqnarray*}
b\1 x S(b\2) & = & \lambda^{-1} b\1w S(b\2) w^{-1} E_{M_1}(e_2 x S(b\3)) \\
& = & \lambda^{-1} b\1 S(b\2) E_{M_1}(e_2 x S(w^{-1} b\3 w)) \\
& = & \lambda^{-1} \eps_t(b\1) E_{M_1}(S(w^{-1})b\2 x e_2)w \\
& = & \lambda^{-1} E_{M_1}(S(w^{-1}) b x e_2) w. 
\end{eqnarray*}
Next note that $\Delta(v') = (1 \otimes v')\Delta(1)$
for all $ v' \in W$,  
which follows from  an application of $S$ to Lemma (\ref{delta on V}). 
Then  let $b' = S(w^{-1})b$ and compute:
$$
b' \lact x = (S(w)b')\1 x S((S(w)b')\2) w^{-1} = 
b'\1 x S(S(w)b'\2) w^{-1} = b'\1 x S(b'\2). \qed 
$$
\renewcommand{\qed}{}\end{proof}

\begin{theorem}
\label{psi theorem}
The mapping $\psi:  x \# b \mapsto xb \in M_2$ defines an isomorphism between 
the algebra $M_2$ and the smash product algebra $M_1 \# B$. 
\end{theorem}
\begin{proof}
That $\psi$ is a linear isomorphism follows from Lemma (\ref{M_2 is M_1B}). 

That $\psi$ is a homomorphism follows almost directly from Eq.\
(\ref{WHA theory}) and 
the conjugation formula in Proposition (\ref{simplification}):
$$ bx = b\1 x \eps_s(b\2) = (b\1 \lact x)b\2, 
$$
since for all $b' \in B$:  $\eps_s(b') = S(b'\1)b'_2 \in W = C_{M_2}(M_1)$.
\end{proof} 

\subsection*{Action of $\mathbf{A}$ on $\mathbf{M}$}
In this subsection, we define a left action of $A$ on $M$ by a formula
similar to that for $\lact$ of $B$ in Proposition (\ref{simplification}).
Denote the antipode of $A$ by $S$ below. Let $\eps_s$ and $\eps_t$
again denote the source and target counital maps on $A$.  

\begin{lemma} 
\label{WHA fact} 
The map $\eps_t$ is a  non-unital module homomorphism 
${}_A A \to {}_{\rm ad}A$ with respect to the left regular  
and adjoint actions of $A$ on itself: for all $a,a' \in A$, 
$$a\1 \eps_t(a') S(a\2) = \eps_t(aa').$$
\end{lemma}

The proof of this and a similar fact for $\eps_s: A_A \to A_{\rm ad}$ is easy
and omitted.
\begin{proposition}
\label{A acts on M}
The mapping $\lact:\, A \otimes M \to M$ given by
\begin{equation}
a \lact m = a\1 m S(a\2)
\end{equation}
is a weak Hopf algebra action of $A$ on $M$. 
\end{proposition}
\begin{proof}
First we check that $a \lact m \in M$ given $m \in M, a \in A$. 
Let $\rho: M_1 \rightarrow M_1 \otimes A$, $\rho(x) = x^{(0)} \otimes x^{(1)}$,
denote the coaction dual
to the action $B \otimes M_1 \to M_1$ above.  
Then $b \lact x = x^{(0)} \la x^{(1)}, b \ra$.
It follows from Eq.\ (\ref{standard action}) that $\rho$ restricted to $A$ is the
comultiplication:  
$$ a^{(0)} \otimes a^{(1)} = a\1 \otimes a\2. $$
Since $M$ is shown above to be  the invariant subalgebra of this action of $B$ on $M_1$, it
is also precisely the coinvariant subalgebra of $\rho$.  
We then compute using Lemma (\ref{WHA fact}): 
\begin{eqnarray*}
\rho(a \lact m)  & = & a\1 m^{(0)} S(a\4) \otimes a\2 \eps_t(m^{(1)}) S(a\3) \\
                & = & a\1 m^{(0)} S(a\3) \otimes \eps_t(a\2 m^{(1)})  \\
                & = & (a \lact m) ^{(0)}\otimes \eps_t((a \lact m)^{(1)}),
\end{eqnarray*}
whence $a \lact m \in M$. 

Since $\eps_s(A) = V = C_{M_1}(M)$, we compute that $\lact$ measures $M$:
\begin{eqnarray*}
(a\1 \lact m)(a\2 \lact m') &=& a\1 m S(a\2)a\3 m' S(a\4) \\
                            & = & a\1 \eps_s(a\2) mm' S(a\3) \\
                            & = &  a \lact (mm').
\end{eqnarray*}
 
We note also that $a \lact 1 = \eps_t(a)$ and that
$$ a \lact (a' \lact m) = (aa') \lact m $$
by the homomorphism and anti-homomorphism properties of $\Delta$ and $S$. 
Finally, $1 \lact m = m$ since both $1\1$ and $S(1\2)$ belong to
$V$, while $1\1 S(1\2) = 1_A$. 
\end{proof} 
\begin{theorem}
\label{phi theorem}
The mapping $\phi:  m \# a \mapsto ma \in M_1$ defines an isomorphism between the
algebra $M_1$ and the smash product algebra $M \# A$. 
\end{theorem}
\begin{proof}
That $\phi$ is a linear isomorphism follows from Lemma (\ref{M_2 is M_1B}). 

That $\phi$ is a homomorphism follows from 
the conjugation formula in Proposition (\ref{A acts on M}):
$$ am = a\1 m \eps_s(a\2) = (a\1 \lact m)a\2, 
$$
since for all $a' \in A$:  $\eps_s(a') = S(a'\1)a'\2 \in V = C_{M_1}(M)$.
\end{proof}

\begin{proposition}
\label{fixed subalg}
For the action of $A$ on $M$, we have $N = M^A$.
\end{proposition}
\begin{proof}
If $n \in N$, then for every $a \in A$:
$$a \lact n = a\1 n S(a\2) = \eps_t(a) (1 \lact n) 
= 1\1 \eps_t(a) n S(1\2) = \eps_t(a) \lact n,  $$
using the definition of a module algebra over a weak Hopf algebra.

We similarly compute for each $x \in M^A, a \in A$:
\begin{eqnarray*}
xS(a) & = & \eps_s(a\1) x S(a\2) \\
      & = & S(a\1) (a\2 \lact x) \\
       & = & S(a\1) (\eps_t(a\2) \lact x) \\
      & = &  S(a\1)\eps_t(a\2) 1\1 x S(1\2) = S(a) (1 \lact x) = S(a)x
\end{eqnarray*}
From the bijectivity of $S: A \to A$ and $e_1 \in A$,
it follows that $e_1 x = x e_1$, so that $xe_1 = e_1 x e_1 = E(x) e_1$,
whence $x = E(x) \in N$. 
\end{proof} 
\end{section}
\begin{section}{Appendix:  The Composite basic construction and a Depth 
$2$ example} 
In this appendix we discuss the two unrelated topics in the title.

Extending the Jones tower in (\ref{Jones tower}) indefinitely
to the right via iteration of the
basic construction for a  subfactor $N \subseteq M$ of positive
index $\lambda^{-1}$, 
Pimsner and Popa \cite{PP2} have shown that the basic construction
of the composite conditional expectation
$$F_n := E \circ E_M \circ \ldots \circ E_{M_{n-1}}: \ M_n \to N$$
is isomorphic to $M_{2n+1}$ with Jones idempotent $f_n \in M_{2n+1}$
given by 
\begin{equation}
\label{PP idempotent}
 f_n = \lambda^{-n(n+1)/2} (e_{n+1} e_n \cdots e_1)(e_{n+2} e_{n+1} \cdots
e_2)\cdots(e_{2n+1} e_{2n} \cdots e_{n+1}). 
\end{equation}
We will prove here that the same is true in the more general algebraic
situation where $M/N$ is a strongly separable extension of index
$\lambda^{-1}$. We do not
need a Markov trace here. This appendix is not needed in Sections 3 and 4. 

Let $F_{M_n} = E_{M_n} \circ \cdots \circ E_{M_{2n}}:
M_{2n+1} \to M_n$. 

\begin{proposition}
The element 
$f_n$ is an idempotent satisfying the characterizing properties
of the basic construction:
$$M_{2n+1} = M_n f_n M_n, $$
$$f_n x f_n = f_n F_n(x) = F_n(x) f_n, \ \ \forall x \in M_n,$$
$$ F_{M_n}(f_n) = \lambda^{n+1} 1. $$
\end{proposition}
\begin{proof}
The proof in \cite{PP2} that $f_n^2 = f_n$, $ F_{M_n}(f_n) = \lambda^{n+1} 1 $
and $ f_n F_n(x) = F_n(x) f_n$ 
is valid here as it only makes use of the $e_i$-algebra $A_{n, \lambda}$,
the subalgebra of $M_n$
 $k$-generated by $e_1, \ldots, e_n$, and an obvious
involution on it.  Note that the theorem is true for $n = 0$
(where $f_0 = e_1$).  Assume inductively that the proposition holds
for $n - 1$ and less.  We use the induction hypothesis 
in the second step below, and the Pimsner-Popa identities 
for sets $f_{n-1}M_{2n-1} =
f_{n-1}M_{n-1}$ in the fifth step:  
\begin{eqnarray*}
 M_{2n+1}  & = & M_{2n} e_{2n+1} M_{2n} \\
 & = & M_{2n-1} e_{2n} M_{2n-1} e_{2n+1}M_{2n-1} e_{2n} M_{2n-1} \\
& = & M_{2n-1} e_{2n}e_{2n+1} M_{n-1} f_{n-1} M_{n-1}e_{2n}M_{2n-1}  \\
& = & M_{2n-2} e_{2n-1}e_{2n} e_{2n+1} M_{2n-2} f_{n-1} M_{2n-2} e_{2n}e_{2n-1} M_{2n-2} \\
& = & M_{2n-2} e_{2n-1}e_{2n} e_{2n+1} f_{n-1}e_{2n}e_{2n-1} M_{2n-2} \\
& = & \cdots =    M_n e_{n+1} \cdots e_{2n+1} f_{n-1} e_{2n} \cdots e_{n+1} M_n = M_n f_n M_n,
\end{eqnarray*}
where the last step is by \cite[Lemma 2.3]{PP2}. 

Let $\tau^{2}$ denote the shift map of $A_{n,\, \lambda} \to 
A_{n+2,\, \lambda}$ induced by $e_i \mapsto e_{i+2}$.  
It follows from the induction hypothesis that $\tau^{2}(f_{n-1})$
is the Jones idempotent for the composite expectation 
$$\widehat{F_{n-1}} := E_{M_1} \circ \cdots \circ E_{M_n}: M_{n+1} \to M_1.$$
Let $x \in M_n$ and $x' = E_{M_{n-1}}(x)$. For the
computation below, we note that $e_{n+1} x e_{n+1} = x' e_{n+1}$
and by \cite[Remark 2.4]{PP2}:
$$ f_n = \lambda^{-n} (e_{n+1} e_n \cdots e_1) \tau^2(f_{n-1})(e_2 e_3 \cdots e_{n+1}). $$
We compute:
\begin{eqnarray*}
\lefteqn{f_n x f_n =}\\
& = &  \lambda^{-2n} (e_{n+1}  \cdots e_1)\tau^2(f_{n-1})
 (e_2 \cdots e_{n+1})x'(e_{n+1}  \cdots e_1) \tau^2(f_{n-1})(e_2  \cdots e_{n+1}) \\
& = & \lambda^{-2n} (e_{n+1}  \cdots e_1) \widehat{F_{n-1}}(e_2  \cdots e_n x' e_{n+1} e_n \cdots e_2 e_1) \tau^2(f_{n-1})(e_2  \cdots e_{n+1}) \\
& = & \lambda^{-n} (e_{n+1}  \cdots e_1)E_M \circ \cdots \circ E_{M_{n-1}}(x)e_1 \tau^2(f_{n-1})(e_2 e_3 \cdots e_{n+1}) \\
& = & F_n(x) f_n. \qed
\end{eqnarray*}
\renewcommand{\qed}{}\end{proof} 

\begin{remark}
It was shown in \cite{NV2} that if $N\subseteq M$ is a II$_1$ subfactor
of finite index and {\em arbitrary} finite depth (see \cite{GHJ} for a definition)
then there exists $k\geq 0$ such that for all $i\geq k$ subfactors 
$N\subseteq M_i$ have depth $2$. It would be interesting to extend this
property to the purely algebraic case (the finite depth property
in this setting was defined in \cite{KN}).
\end{remark}

As a final topic in this appendix we provide examples of depth $2$ extensions
 in the next  proposition and corollary.

\begin{proposition}
Suppose $M/N$ is a weakly irreducible, symmetric,
strongly separable extension
such that its bimodule projection $E: M \to N$ has dual
bases in the centralizer $U$.  Suppose
moreover that the center $C$ of $U$ coincides with
the center $Z$ of $N$. Then $M/N$ has depth $2$.  
\end{proposition}
\begin{proof}
Let $x_i, y_i \in U = C_M(N)$ be dual bases of $E$.  
It follows that $M \cong N \otimes_Z U$ via $m \mapsto E(mx_i) \otimes y_i$.
By the symmetry condition on $E$, $E$ restricted to $U$ is a trace
with values in $Z = C$.
Then $\lambda x_i \otimes y_i$ is the symmetric separability element
and $$ u \mapsto \lambda x_i u y_i $$
gives a $C$-linear projection of $U$ onto $C$ coinciding with
$ E|_U$, since
$U$ is an Azumaya $C$-algebra  
\cite[Section 3]{SK}.

Let $z_i = \lambda^{-1} x_i e_1$ and $w_i = e_1 y_i$ in $M_1$:
these are dual bases of $E_M: M_1 \to M$ by the Basic Construction Theorem.
But we see that $z_i , w_i \in A$.

Next we compute that there are dual bases $x'_i, y'_i \in V = C_{M_1}(M)$
 for $E_M$.  By the construction of the last paragraph, it
follows that $E_{M_1}$ has dual bases in $B$, whence
$M/N$ has depth $2$.
We let $x'_i = x_j x_i e_1 y_j$ and $y'_i = x_k y_i e_1 y_k$,
both in $V$.
It suffices to compute for $a, b \in M$:
\begin{eqnarray*}
E_M(ae_1b x'_i)y'_i & = & E_M(aE(bx_jx_i)e_1 y_j) y'_i \\
                    & = & \lambda a E(x_i b x_j)y_j x_k y_i e_1 y_k \\
                    & = & \lambda a x_i x_k y_i e_1 y_k b  \\
                    & = & a  e_1 E(x_k) y_k b = ae_1 b.
\end{eqnarray*}
Similarly we compute $x'_i E(y'_i ae_1 b) = ae_1 b$ by using the equivalent
expressions $x'_i = x_j e_1 x_i  y_j$ and $y'_i = x_k e_1 y_i  y_k$.
\end{proof}
For the next corollary-example, we need a few definitions.  
An algebra $A$ is \textit{central} if its center is trivial, $Z(A) = k1$. 
A ring extension $M/N$ is \textit{H-separable} (after Hirata) if there are
elements $f_i \in (M \otimes M)^N$ and $u_i \in U = C_M(N)$ such that 
$e_1 = u_i f_i$,  where $e_1$ again denotes $1 \otimes 1$ in $M \otimes_N M$
\cite{K2}.
\begin{corollary}
Suppose $M/N$ is a split H-separable extension of central algebras
where $U$ is Kanzaki separable. Then $M/N$ is a depth
$2$ strongly separable extension.  
\end{corollary}
\begin{proof}
By the results of \cite[Theorem 2.1]{XY}, 
the center of $U$ is trivial and $ N \otimes U \cong M$ via $n \otimes u \mapsto
nu$ for $n \in N, u \in U$. But by the hypothesis $U$ has  non-degenerate trace
$t: U \rightarrow k$ with dual bases $x_i, y_i \in U$.
It follows that $E: M \to N$ defined by $E(nu) = \lambda nt(u)$,
where $\lambda^{-1} = t(1)$, 
has dual bases in $U$.  The conclusion now follows readily
from the proposition. 
\end{proof}

\end{section}

\bibliographystyle{amsalpha}

\begin{thebibliography}{A}

\bibitem{BM}
R.~Blattner and S.~Montgomery,
\textit{A duality theorem for Hopf module algebras}, 
\textit{J. Algebra} \textbf{95} (1985), 153--172.

\bibitem{BNSz} G.~B\"ohm, F.~Nill and K.~Szlach\'anyi,
Weak Hopf algebras, I. Integral theory and $C^*$-structure,
\textit{J. Algebra} \textbf{221} (1999), 385-438, doi:10.1006/jabr.1999.7984. 

\bibitem{BSz} G.~B\"ohm and  K.~Szlach\'anyi,   
A coassociative $C^*$-quantum group with nonintegral dimensions,
\textit{Lett. in Math. Phys}, \textbf{35}  (1996), 437--456.


\bibitem{EN} P.~Etingof and  D.~Nikshych,
Dynamical quantum groups at roots of $1$, 
to appear in \textit{Duke Math. J.},
\texttt{math.QA/0003221} (2000).

\bibitem{GHJ} F.~Goodman, P.~de la Harpe, and V.F.R.~Jones, 
``Coxeter Graphs and Towers of Algebras,''  M.S.R.I. Publ.\  
\textbf{14}, Springer, Heidelberg, 1989.

\bibitem{HS} K.~Hirata and K.~Sugano, 
On semisimple extensions and separable extensions over 
non commutative rings,  \textit{J. Math.\ Soc.\ Japan}  
\textbf{18} (1966), 360--373.

\bibitem{J1} V.F.R.~Jones,
Index for subfactors, \textit{Inventiones Math.} \textbf{72} (1983), 1--25.

\bibitem{J2} V.F.R.~Jones,
Index for subrings of rings, {\it Contemp.\ Math.}{ \bf 43} A.M.S., 
(1985), 181--190.

\bibitem{K1} L.~Kadison, 
The Jones polynomial and certain separable Frobenius extensions, 
\textit{J. Algebra} \textbf{186} (1996), 461--475, doi: 10.1006/jabr.1996.0383.  

\bibitem{K2} L.~Kadison, 
``New examples of Frobenius extensions,'' 
University Lecture Series \textbf{14},
Amer.\ Math.\ Soc., Providence,  1999. 

\bibitem{KN} L.~Kadison and D.~Nikshych,
Outer actions of  centralizer Hopf algebras on separable extensions,
 \textit{Comm.\ Alg.}, to appear.

\bibitem{Kan} T.~Kanzaki,
Special type of separable algebra over commutative ring,
{\it Proc.\ Japan Acad.} {\bf 40} (1964), 781--786.

\bibitem{Kasch1} F.~Kasch,
Projektive Frobenius Erweiterungen,  \textit{Sitzungsber.\ Heidelberg.\
 Akad.\ Wiss.\ Math.-Natur.\ Kl.} (1960/1961), 89--109.

\bibitem{Kasch2}F.~Kasch, 
Dualit\"{a}tseigenschaften von Frobenius-Erweiterungen,  
{\it Math.\ Zeit.} {\bf 77 } (1961), 219--227. 

\bibitem{M} S.~Montgomery,
``Hopf algebras and their actions on rings,''  CBMS Regional Conf.\ 
Series in Math.\ \textbf{82}, A.M.S., Providence, 1993. 

\bibitem{N} D.~Nikshych, 
A duality theorem for quantum groupoids, in: ``New Trends in Hopf
Algebra Theory,'' eds.\ N. Andruskiewitsch, F. Santos and H.-J. Schneider,
\textit{Contemp.\ Math.} \textbf{267}
(2000), 237--243. 
 
\bibitem{NTV} D.~Nikshych, V.~Turaev, and L.~Vainerman,
Quantum groupoids and invariants of knots and 3-manifolds, preprint,
\texttt{math.QA/0006078} (2000).

\bibitem{NV1} D.~Nikshych and L.~Vainerman,
A characterization of depth $2$ subfactors of II${}_1$  factors,
\textit{J. Func. Analysis} \textbf{171} (2000), 278--307.

\bibitem{NV2} D.~Nikshych and  L.~Vainerman,
A Galois correspondence for actions of quantum  groupoids on II$_1$-factors,
\textit{J. Func. Analysis}, \textbf{178} (2000), 113-142.

\bibitem{NV3} D.~Nikshych and L.~Vainerman,
Finite dimensional quantum groupoids
and their applications, to appear in the Proceedings of 
``Hopf Algebras'' workshop, MSRI Publications (2000), 
\texttt{math.QA/0006057}.

\bibitem{Oc} A.~Ocneanu,
\textit{Quantized groups, string algebras and Galois theory for algebras},
Operator Algebras and Applications, Vol. 2, 
London Math. Soc. Lecture Notes Series \textbf{135}, Cambridge Univ.
Press, Cambridge, U.K., (1988).

\bibitem{O} T.~Onodera,
Some studies on projective Frobenius extensions,
{\it J. Fac.\ Sci.\ Hokkaido Univ.\ Ser.\ I}, {\bf 18} (1964), 89-107.

\bibitem{PP2} M.~Pimsner and S.~Popa,
Iterating the basic construction,
\textit{Trans.\ AMS} \textbf{310} (1988), 127--133.

\bibitem{SK} A.A.~Stolin and L.~Kadison,
Separability and Hopf algebras, in: ``Algebra and its Applications,''
eds.\ Huynh, Jain, Lopez-Permouth, Contemp.\ Math.\ vol.\ \textbf{259}, AMS,
Providence, 2000, 279--298.

\bibitem{Sz} K.~Szlach\'anyi,
Weak Hopf algebra symmetries of $C^*$-algebra inclusions, preprint
\texttt{math.QA/0101005} (2001).

\bibitem{S} W.~Szyma\'nski,
Finite index subfactors and Hopf algebra crossed products,
\textit{Proc.\ Amer.\ Math.\ Soc.} \textbf{120} (1994), no. 2, 519--528.


\bibitem{W} Y.~Watatani,
Index of $C^{\ast}$-subalgebras,
{\it Memoirs A.M.S.} {\bf 83} (1990).

\bibitem{XY} J.~Xiaolong and X.~Yongchua,
H-separable rings and their Hopf-Galois extensions,
\textit{Chin.\ Ann.\ Math.} \textit{19B} (1998), 311-320.

\bibitem{Y} K.~Yamagata,
Frobenius Algebras, in: {\it Handbook of Algebra, Vol.\ 1}, ed.\ M. Hazewinkel,
Elsevier, Amsterdam, 1996, 841--887.

\end{thebibliography}

\end{document}